\documentclass[11pt,a4paper]{paper}
\usepackage{fullpage}
\usepackage[english]{babel}
\usepackage{amsmath}
\usepackage{amssymb}
\usepackage{graphicx}
\usepackage{tikz,pgfplots,pgfplotstable}
\usepackage{color}
\usepackage{hyperref}             
\hypersetup{linktocpage}

\newtheorem{de}{Definition}[section]
\newtheorem{theo}{Theorem}[section]
\newtheorem{cor}[theo]{Corollary}
\newtheorem{prop}[theo]{Proposition}
\newtheorem{lem}[theo]{Lemma}
\newtheorem{rem}[theo]{Remark}

\numberwithin{figure}{section}

\newcommand{\Proof}{\medskip\noindent {\bf Proof. }}

\newcommand{\A}{{\cal A}}
\newcommand{\B}{{\cal B}}
\newcommand{\D}{{\cal D}}
\newcommand{\C}{{\cal C}}
\newcommand{\E}{{\cal E}}

\newcommand{\cH}{{\cal H}}

\newcommand{\N}{\mathbb{N}}
\newcommand{\cO}{{\cal O}}
\newcommand{\Ra }{{\cal R}}
\newcommand{\R}{\mathbb{R}}
\newcommand{\Sa}{{\cal S}}

\newcommand{\al}{\alpha}
\newcommand{\be}{\beta}
\newcommand{\ga}{\gamma}

\newcommand{\e}{\varepsilon}
\newcommand{\Om}{\Omega}
\newcommand{\om}{\omega}
\newcommand{\ph}{\varphi}
\newcommand{\Ph}{\Phi}
\newcommand{\s}{\sigma}
\newcommand{\Si}{\Sigma}

\newcommand{\htwo}{{H^{2}(\Om)}}
\newcommand{\hh}{{H^{2}(\Om)}}

\newcommand{\BV}{\text{BV}(\Om)}
\newcommand{\ltwo}{{L^{2}(\Om)}}

\newcommand{\linf}{{L^{\ii}(\Om)}}
\newcommand{\hpartial}{{H^{1/2}(\dr\Om)}}

\newcommand{\1}{{\bf{1}}}
\newcommand{\ii}{\infty}
\newcommand{\wh}{\widehat}
\newcommand{\ton}{\text{ \ on }}
\newcommand{\tin}{\text{ \ in }}
\newcommand{\dist}{\text{dist}}
\newcommand{\nn}{\nonumber}
\newcommand{\bs}{\backslash}
\newcommand{\dr}{\partial}

\newcommand{\g}{\nabla}
\newcommand{\p}{\cdot}
\newcommand{\ti}{\tilde}
\newcommand{\ul}{\underline}
\newcommand{\ol}{\overline}
\newcommand{\su}{\subset}
\newcommand{\La}{\triangle}
\newcommand{\supp}{\text{supp\ }}
\newcommand{\mat}{\text{mat}}

\newcommand{\ra}{\longrightarrow}
\newcommand{\mt}{\longmapsto}
\newcommand{\cqfd}{\hfill $\square$\\ \medskip}
\newcommand{\n}[2]{\left\|{#1}\right\|_{#2}}
\newcommand{\io}{\int_\Om}

\newcommand{\bb}{\hspace{1 cm}}

\title{Reconstruction and stability in acousto-optic imaging for absorption maps with bounded variation\thanks{\footnotesize This work was supported  by the
ERC Advanced Grant Project MULTIMOD--267184.}} 

\author{Habib Ammari\thanks{\footnotesize Department of Mathematics and Applications,
Ecole Normale Sup\'erieure, 45 Rue d'Ulm, 75005 Paris, France
(habib.ammari@ens.fr, laurent.seppecher@ens.fr).} \and
Loc Hoang Nguyen\thanks{Mathematics Section, \'Ecole Polytechnique F\'ed\'erale de Lausanne,
Station 8,  CH-1015 Lausanne, Switzerland (loc.nguyen@epfl.ch).}  \and Laurent Seppecher\footnotemark[2] }

\begin{document}
\maketitle

\begin{abstract}
The aim of this paper is to propose for the first time a reconstruction scheme and a stability result for recovering from acoustic-optic data absorption distributions with bounded variation.
The paper extends earlier results in \cite{chap_2} and \cite{chap_3} on smooth absorption distributions. It opens a door for a mathematical and numerical framework for imaging,  from internal data,  parameter distributions with high contrast in biological tissues. 
\end{abstract}

\bigskip

\noindent {\footnotesize Mathematics Subject Classification
(MSC2000): 35R30, 35B30.}


\section{Introduction}
In the recent papers \cite{chap_2},\cite{chap_1}, and \cite{chap_3},  an original mathematical and numerical framework for modeling biomedical imaging modalities based on mechanical perturbations of the  medium is developed. The objective is to enhance the resolution and stability of tissue property imaging. 

Many kinds of waves propagate in biological tissues over certain frequency ranges. Each one of them can be used to provide an image of a specific physical parameter. Low-frequency electromagnetic waves are sensitive to electrical conductivity; optical waves tell about optical absorption, ultrasonic waves reveal tissue's density, mechanical shear waves indicate how tissues respond to shear forces. However, single-wave imaging modalities are known to suffer from low specificity as well as intrinsic instabilities and low resolution; see \cite{AMMARI-08} and \cite{seobook}. These fundamental deficiencies are impossible to eliminate, unless additional a priori information is incorporated. Single-wave imaging modalities can only be used for anomaly detection.  Expansions techniques for data analysis, which reduce the set of admissible solutions and the number of unknowns, allow robust and accurate reconstruction of the location and of some geometric features of the anomalies, even with moderately noisy data.

One promising way to overcome the inherent limits of single-wave imaging and provide a stable and quantitative reconstruction of a distribution of  physical parameters is to combine different wave-imaging modalities; see again \cite{AMMARI-08} and \cite{seobook}. A variety of multi-wave imaging approaches are being introduced and studied. In such approaches, two or more types of physical waves are involved in order to overcome the individual deficiencies of each one of them and to combine their strengths.  Because of the 
way the waves are combined, multi-wave imaging can produce a single image with the best contrast and resolution properties of the two waves. 

Three different types of wave interaction can be exploited in multi-wave imaging  \cite{fink}: (i) the interaction of one kind of wave with tissue can generate a second kind of wave; (ii) a low-frequency wave that carries information about the desired contrast can be locally modulated by a second wave that has better spatial resolution; (iii) a fast propagating wave can be used to acquire a spatio-temporal sequence of the propagation of a slower transient wave. 

In \cite{chap_2} and \cite{chap_3},  by  mechanically perturbing the medium we proved both analytically and numerically  the stability and resolution enhancement for reconstructing optical tissue parameters. We showed how the high contrast of  optical tomography \cite{simon} can be coupled to the high resolution of the acoustic propagation in soft tissues. The use of mechanical perturbations of the medium modeled by acoustics equations in fluids enhance the resolution to the order of the front width of the acoustic wave, which propagates inside the object. It dramatically increases the low resolution of optical tomography \cite{john}.

This paper is a continuation and an extension of the work started in \cite{chap_2} and \cite{chap_3}. We keep here the same models for the diffusive light propagation \cite{born} and for the acoustic perturbations. Our aim is to extend the reconstruction algorithm developed in  \cite{chap_2} to a large class of non smooth functions taken in a subclass of $\BV$, the set of functions with bounded variation. 

The reconstruction and the stability of the inversion are shown in this general case. Such an extension is essential for applying the proposed hybrid method to biological tissues. Indeed, the physiologic parameters that we want to recover cannot be considered smooth or piecewise smooth as assumed in \cite{chap_3}.

Under this natural assumption, new mathematical difficulties rise to prove that the acousto-optic data contain enough information for reconstructing the absorption map. The lack of smoothness also causes difficulties to ensure the stability of the algorithm. This paper resolves these challenging issues. It provides both an original reconstruction formula and a new stability result in the general setting. As far as we know, together with the recent work \cite{otmar}, it is the first work in imaging discontinuous parameter distributions from internal measurements.

Throughout this paper, we denote by $\mathcal{S}$ the space of Schwartz and by $\mathcal{S}^\prime$ its dual.  We use the notation $H^s$ for the usual Sobolev spaces and set $\D$ to be the set of $\mathcal{C}^\infty$ compactly supported functions. 

As in \cite{chap_2} and \cite{chap_3}, we consider a smooth bounded domain $\Om$ of $\R^d$, for $d\in\{2,3\}$, and a light fluence field defined as the unique solution of the diffusion equation

\begin{equation}\label{sys_Phi}\left\{\begin{aligned}
-\La\Ph+a\Ph &=0\ \ \ \tin\Om,\\
l\dr_\nu\Ph+\Ph &=g\ \ \ton\dr\Om,
\end{aligned}\right.\end{equation}
where $a\in\linf$ satisfying $a\geq\ul a>0$ and $\supp(a-a_0)\su D \Subset \Om$ is the absorption
 parameter to be recovered; see \cite{simon} and \cite{john}. The extrapolation length $l$, and the bounds $\ul a$ and $a_0$ are known positive constants. The incoming illumination $g\in\hpartial$ is a non negative non zero map and is also supposed to be known. Moreover, the support $D$ of $a-a_0$ is assumed to be smooth. 

 The acoustic perturbations are assumed to be generated by spherical pressure waves. Let $\eta$ be the front width of the acoustic wave and let $w$ be the wave shape. The acoustic perturbations take the form:

\begin{equation}\label{eq_defv}\begin{aligned}
v_{y,r,\eta}(x)=\frac{\eta}{r}w\left(\frac{|x-y|-r}{\eta}\right)\frac{x-y}{|x-y|},\bb\forall \; x\in\R^d\bs\{y\},
\end{aligned}\end{equation}
where $y\in Y\su\R^d$, $\eta>0$ and $r\in]\eta,+\ii[$; see  \cite{chap_1}. 
Here, $Y \subset \Omega \setminus \overline{D}$ is a smooth surface. Moreover,  the map $w\in\D(\R)$ is non negative and satisfies $\supp(w)\su[-1,1]$, $w'>-1$ and $\| w\|_{L^1}=1$. The last assumption ensures that the map $x\mt x+v(x)$ is a diffeomorphism. 

The effect of the displacement $v$ on the absorption map is assumed to be only a shifting effect, that is, to say that $a$ becomes $a_v$ implicitly defined on $\Om_v=(Id+v)(\Om)$ by

\begin{equation}\begin{aligned}
a_v(x+v(x))=a(x),\bb\forall \; x\in\Om_v,
\end{aligned}\end{equation}
or equivalently, by the formula $a_v=a\circ(Id+v)^{-1}$. We introduce the displaced light fluence as the unique solution of 

\begin{equation}\label{sys_Phimod}\left\{\begin{aligned}
-\La\Ph_v+a_v\Ph_v &=0\ \ \ \tin\Om,\\
l\dr_\nu\Ph_v+\Ph_v &=g\ \ \ton\dr\Om,
\end{aligned}\right.\end{equation}
by extending $a_v$ by $a_0$ if necessary. Computing now the cross-correlation on the boundary $\dr\Om$ between $\Ph$ and $\Ph_v$ it follows that

\begin{equation}\label{eq_corr4}\begin{aligned}
\frac 1 l\int_{\dr\Om}(\Ph-\Ph_v)g=\int_\Om(a_v-a)\Ph\Ph_v.
\end{aligned}\end{equation}
Assume that the term in the left-hand side of the above identity can be measured. We define the measurement as the real quantity given by
\begin{equation}\label{}\begin{aligned}
M_v=\frac 1{\eta^2}\int_\Om(a_v-a)\Ph\Ph_v. 
\end{aligned}\end{equation}
Throughout this paper, we assume that $M_v$ is known for any displacement field $v$ given by (\ref{eq_defv}). 

For a smooth surface $Y\su \Om \setminus \overline{D}$ and $\eta>0$, we assume that we are in possession of  

\begin{equation}\begin{aligned} \label{imagdata}
M_\eta(y,r)=\frac 1{\eta^2}\int_\Om(a_{v_{y,r,\eta}}-a)\Ph\Ph_{v_{y,r,\eta}},\bb\forall \; (y,r)\in Y\times]\eta,+\ii[ .
\end{aligned}\end{equation}

The imaging problem considered in this paper is to reconstruct $a$ from the measurement data $M_\eta$ given by (\ref{imagdata}). The aim is to prove that the reconstruction algorithm from acousto-optic differential measurements presented in \cite{chap_2} can be extended for a very general class of discontinuous absorption maps.  For doing so, we start  from the same differential boundary measurements (\ref{imagdata}) and consider the case where $a$ has bounded variations. Under some additional hypothesis, we correctly interpret the first order term in the asymptotic formula when $\n{v}{L^\ii}$ goes to zero.  Then, by giving a weak definition of the spherical means Radon transform $\Ra$, we show how the internal data $\Psi$, satisfying
\begin{equation}\nn\begin{aligned}
\Ph^2\text{D} a=\text{D} \Psi+\g\times G,
\end{aligned}\end{equation}
can be reconstructed stably in $H^s(D)$ with $s<1/2$ and $D$ being a smooth domain.
This is done through a stable reconstruction of $\Ra[\Psi]$ in $H^{(d-1)/2+s}$. Here, $\text{D} a$ and $\text{D} \Psi$ are defined by (\ref{defda}).

%


The second part is to show that a stable reconstruction of the absorption map $a$ is possible from this internal data $\Psi$. In order to do so, we establish a system of two coupled elliptic equations for $(a,\Phi)$ and solve this coupled system by the classical fixed point theorem. We also show that the solution depends continuously on $\Psi$ and therefore can verify the global stability of the reconstruction. 

Finally, we present numerical illustrations to substantiate the potential of the proposed method. 
We consider the imaging of a highly discontinuous absorption map, chosen from a real biological tissue data.

\section{Preliminaries}

In order to work with a wide set of discontinuous functions, we introduce $BV(\Om)$ and several important subspaces of $BV(\Om)$.
\bigskip

\subsection{Some subclasses of functions with bounded variation}
\begin{de} A function $u\in L^1(\Om)$ is said to have bounded \index{bounded variation} variation if its weak derivative $D u$ is a finite \index{Radon measure} Radon measure. For any $\ph\in\C^1_c(\Om)^d$, we have
\begin{equation}\nn\begin{aligned}
\int_\Om u(x)\g\cdot \ph(x)dx=-\int_\Om \ph(x)\cdot D u(dx).
\end{aligned}\end{equation} 
\end{de} 
The Radon measure $Du$ can be uniquely decomposed into three singular measures as follows:
\begin{equation} \begin{aligned} \label{defda}
D u=D _lu+D _ju+D _cu,
\end{aligned}\end{equation}
which are respectively called the Lebesgue part, the jump part, and the Cantor part of $D u$. The Lebesgue part is absolutely continuous with respect to the Lebesgue measure and is identified to $D_lu\in L^1(\Om)^d$, which is called the smooth variation of $u$. The jump part $D _ju$ is such that there exists a set $S\su\Om$ of Hausdorff dimension $(d-1)$, rectifiable admitting the existence of a generalized normal vector $\nu_S(x)$ for almost every $x\in S$. This part is written as $$D _j u=[u]_S\nu_S\cdot {\cal H}^{d-1}_S,$$ where $[u]_S\in L^1(S,{\cal H}^{d-1}_S)$ is the jump of $u$ over $S$ and ${\cal H}^{d-1}_S$ is the Hausdorff measure on $S$. The Cantor part $D _cu$ is supported on a set of Hausdorff dimension less than $(d-1)$, which means that its ${(d-1)}$-Hausdorff-measure is zero; see \cite{alberti_amb}.

In many cases it is very difficult to deal with such a general measure derivative. We introduce the special class of functions of bounded variation $SBV(\Om)$. This class still describes a very large set of discontinuous functions.

\begin{de} A function $u\in\BV$ is in the special class of bounded variation if $D _cu=0$. We denote by $$SBV(\Om)=\left\{u\in\BV,\ D _cu=0\right\}.$$
\end{de}

In some cases, we shall work in some specific $L^p$ framework. Hence, we  use the following spaces.

\begin{de} For any $p\in[1,+\ii]$, we define
 $$SBV^p(\Om)=\left\{u\in SBV(\Om)\cap L^p(\Om),\ D_l u\in L^p(\Om)^d,\ [u]_S\in L^p(S,{\cal H}^{d-1}_S)\right\}.$$
\end{de}

Roughly speaking, a function $u\in SBV^p(\Om)$ is a function of class $W^{1,p}$ admitting surface discontinuities. In the following, we state some \index{Sobolev regularity result} Sobolev regularity results for functions of bounded variation. The embedding rule for  $BV(\Om)$ in the Sobolev spaces behaves like that for $W^{1,1}(\Om)$.

\begin{prop}[$BV(\Om)$ embedding in Sobolev spaces] For any $s\in \R^+, p\geq 1$, if $W^{1,1}(\Om)\hookrightarrow W^{s,p}(\Om)$ continuously, then 
$BV(\Om)\hookrightarrow W^{s,p}(\Om)$ continuously.
\end{prop}

If a function is in $SBV^\ii(\Om)$ we can expect a better Sobolev regularity. We provide the following embedding result.

\begin{prop}\label{prop_sbvinf2sobolev} For any $0\leq \al<\frac 1 2$, $SBV^\ii(\Om)\hookrightarrow H^\al(\Om)$.
\end{prop}
\Proof Consider $u\in SBV^\ii(\Om)$. $D u=D_l u+[u]_S\nu_S {\cal H}^{d-1}_S$ where $S$ is a rectifiable surface, $D_l u\in L^\ii(\Om)^d$ and  $[u]_S\in L^\ii(S,{\cal H}^{d-1}_S)$. We introduce a continuous trace operator $\ga_S:H^{1-\al}(\Om)\ra L^2(S)$ and consider a test function $\ph\in\D(\Om)^d$ to write
\begin{equation}\nn\begin{aligned}
\left<D u,\ph\right>_{\D'(\Om)^d,\D(\Om)^d} &=\int_\Om D_l u\cdot \ph+\int_S [u]_S\nu_S\cdot \ph{\cal H}^{d-1}_S\\
\left|\left<D u,\ph\right>_{\D'(\Om)^d,\D(\Om)^d}\right| &\leq\n {D_l u}{L^\ii(\Om)}\n {\ph}{L^2(\Om)}+\n {[u]_S}{L^\ii(S)}\n {\ph}{L^2(S)}\\
&\leq \n{\gamma_S}{{\cal L}(H^{1-\al}(\Om),\ltwo)}\left(\n {D_l u}{L^\ii(\Om)}+\n {[u]_S}{L^\ii(S)}\right)\n{\ph}{H^{1-\al}(\Om)}.
\end{aligned}\end{equation}
This proves that $D u\in H^{\al-1}(\Om)^d$ and so, $u\in H^\al(\Om)$.\cqfd
\medskip

\subsection{The light fluence operator}\label{subsection 2.2}
The light fluence $\Ph$ associated to the absorption $a$ is defined as the solution of
\begin{equation}\label{eq_light}\left\{\begin{aligned}
-\La\Ph+a\Ph &=0\ \ \ \tin\Om,\\
l\dr_\nu\Ph+\Ph &=g\ \ \ton\dr\Om,
\end{aligned}\right.\end{equation}
where $g$ is smooth (in $H^{3/2}(\partial \Om)$), non negative, and non zero. This problem is well posed if $a\in L^\ii(\Om)$ and admits a positive lower bound. Throughout this  paper, we assume that there exist three constants $0<\ul a\leq a_0\leq\ol a<+\ii$ such that $\ul a\leq a\leq\ol a$ in $\Om$ and $\supp (a-a_0)\su D$. Under this condition, the light fluence $\Phi$ is uniquely determined in $\htwo$. We define the set of the admissible absorption maps by

\begin{equation}\begin{aligned}
\A_0=\left\{a\in\ltwo,\ \ul a\leq a\leq \ol a,\ \supp(a-a_0)\su D\right\}
\end{aligned}\end{equation}
and  the light fluence operator as follows.

\begin{de}\label{def_F} Let the light fluence operator $F$ be given by 
\begin{equation}\nn\begin{aligned}
F:\A_0 &\ra H^2(\Om)\\
a &\mt \Ph,
\end{aligned}\end{equation}
where  $\Ph$ is the unique solution of (\ref{eq_light}).
\end{de}
As in dimensions $2$ and $3$, $\htwo\hookrightarrow\linf$ we define the  following two quantities

\begin{equation}\begin{aligned}
\ul\Ph=\inf_{a\in\A_0}\inf_{x\in\Om} F[a](x),\\
\ol\Ph=\sup_{a\in\A_0}\sup_{x\in\Om} F[a](x).
\end{aligned}\end{equation}

The following result is from \cite{chap_2}.

\begin{prop} The quantity $\ol\Ph$ is finite and depends only on $g$, $l$, $\Om$ and $\ul a$. Moreover, if $g\geq 0$ and $g\neq 0$ in $\dr\Om$, then $\ul\Ph>0$ and depends only on $g$, $l$, $\Om$, and $\ol a$.
\end{prop}

The following proposition is a direct application of standard elliptic regularity results \cite{GilbargTrudinger:1977} on the equation satisfied by $F[a]-F[a']$:
\[
	\left\{\begin{aligned}
-\La (F[a]-F[a'] )+a (F[a]-F[a']) &= (a' - a)F[a']\ \ \ \tin\Om,\\
l\dr_\nu (F[a]-F[a'])+(F[a]-F[a']) &= 0\ \ \ton\dr\Om.
\end{aligned}\right.
\]

\begin{prop}\label{lem_Fcont} The operator $F$ is Lipschitz continuous from $\A_0$ to $\htwo$ in the sense that there exists a constant $C>0$ depending only on $\Om$ such that for any $a$ and $a'$ in $\A$, we have

\begin{equation}\nn\begin{aligned}
\n{F[a]-F[a']}{\hh}\leq C\ol\Ph\n{a'-a}{L^2(\Om)}.
\end{aligned}\end{equation}
\end{prop}

\medskip

In the following, we will suppose that $a$ is in $SBV^\ii(\Om)$ and get from that a little Sobolev regularity enhancement due to Proposition \ref{prop_sbvinf2sobolev}. We have $a\in H^s(\Om)$ for $s\in]0,\frac 1 2[$. For such number $s$, we define a new admissible set for the absorption map:

\begin{equation}\begin{aligned}
\A_s=\left\{a\in\A_0\cap H^s(\Om),\ \n{a}{H^s(\Om)}\leq R_{\A_s}\right\},
\end{aligned}\end{equation}
where $R_{\A_s}$ is a positive real number called the radius of $\A_s$. This gain of regularity for $a$ implies that of regularity for $\Ph=F[a]$, which is stated in the following proposition.

\begin{prop} Assume that $g$ is the trace of a smooth function on $\partial \Omega$. Then for any $s\in]0,\frac 1 2[$ and any $a\in\A_s$, $F[a]\in H^{2+s}(\Om)$. Moreover, the map

\begin{equation}\nn\begin{aligned}
F:\A_s\ra H^{2+s}(\Om)
\end{aligned}\end{equation}
is Lipschitz continuous in the following sense: There exists a constant $C>0$ depending only on $\Om$ and $s$ such that, for any $a$ and $a'$ in $\A_s$, we have

\begin{equation}\begin{aligned}
\n{F[a]-F[a']}{H^{2+s}(\Om)}\leq C  (\ol \Phi + || \nabla \Phi ||_{L^\infty}) \n{a'-a}{H^s(\Om)}.
\end{aligned}
\label{13}\end{equation}
\label{prop 2.5}
\end{prop}

Proposition \ref{prop 2.5} follows immediately from standard regularity estimates. In dimensions $2$ and $3$, $H^2(\Omega) \subset L^{\infty}(\Omega)$. Hence, $\Phi$ satisfies
\[
		\Delta \Phi = a \Phi \in L^{\infty}(\Omega).
\]
This and the smoothness of $g$ imply $\Phi \in \C^{1, \alpha}(\overline{\Om})$ for some $\alpha
 \in (0, 1).$ 
\medskip

\subsection{Spherical means Radon transform} \label{sdfgd}
Here, we introduce  the spherical means Radon transform $\Ra$ and the normalized spherical flow operator $\vec\Ra$. We extend their definition to tempered distributions in order to deal with derivative of non smooth functions. We also give  several useful properties of these operators. We denote by $\Sigma=Y\times ]0,+\ii[$.

\begin{de}[Spherical means Radon transform] For any function $f\in\C^0(\R^d)$, we define its spherical means Radon transform $\Ra[f] \in\C^0(Y\times]0,+\ii[)$ by

\begin{equation}\nn\begin{aligned}
\Ra[f](y,r)=\int_{S^{d-1}}f(y+r\xi)\s(d\xi),\ \ \ \forall \; (y,r)\in\Sigma,
\end{aligned}\end{equation}
where $\s$ is the surface measure of the unit sphere. To extend this definition to distributions, we introduce the dual operator $\Ra^*:\Sa(\Sigma)\ra\Sa(\R^d)$ defined for any $\ph\in\Sa(\Sigma)$ by

\begin{equation}\nn\begin{aligned}
\Ra^*[\ph](x)=\int_Y\frac{\ph(y,|x-y|)}{|x-y|^{d-1}}\s(dy).
\end{aligned}\end{equation}
Then, for any tempered distribution $u\in\Sa'(\R^d)$, we define its spherical mean Radon transform $\Ra[u]\in\Sa'(\Sigma)$ as follows:

\begin{equation}\nn\begin{aligned}
\left<\Ra[u],\ph\right>_{\Sa'(\Sigma),\Sa(\Sigma)}=\left<u,\Ra^*[\ph]\right>_{\Sa'(\R^d),\Sa(\R^d)},\ \ \ \forall \; \ph\in\Sa(\R^d).
\end{aligned}\end{equation}
\end{de}
Injectivity and invertibility issues for $\Ra$ have been studied in several works; see, for instance, \cite{victoire}.  In \cite[Corollary 6.4]{victoire}, the continuity of $\Ra$ and its inverse were proved. The following result holds.

\begin{theo} Consider $s\in\R$ and suppose that for some $\al<s$ and any $u\in H^\al(\Om)$ with compact support, $\Ra[u]=0$ implies $u=0$. Then there exist two positive constants $c_1$ and $c_2$ such that

\begin{equation}\nn\begin{aligned}
\n{u}{H^\al(\Om)}\leq c_1 \n{\Ra[u]}{H^{\al+\frac{d-1}{2}}(\Sigma)}\leq c_2 \n{u}{H^\al(\Om)}.
\end{aligned}\end{equation}
\end{theo}
In the following, we always suppose that we are in the context where this theorem applies. Injectivity issues are essentially controlled by the set of centers $Y$;  
see, for instance, \cite{quinto}.

\begin{de}[Spherical flow operator] For any function $F\in\C^0(\R^d)^d$, we define its normalized flow through the sphere $S(y,r)$,  $\vec\Ra[F]\in\C^0(Y\times]0,+\ii[)$ by

\begin{equation}\label{def_Ravect}\begin{aligned}
\vec\Ra[F](y,r)=\int_{S^{d-1}}F(y+ r \xi)\cdot\xi \s(d\xi),\ \ \ \forall \; (y,r)\in\Sigma.
\end{aligned}\end{equation}

To extend this definition to distributions, we introduce the dual operator $\vec\Ra^*:\Sa(\Sigma)\ra\Sa(\R^d)$ defined for any $\ph\in\Sa(\Sigma)$ by

\begin{equation}\nn\begin{aligned}
\vec\Ra^*[\ph](y,r)=\int_Y\frac{\ph(y,|x-y|)}{|x-y|^d}(x-y)\sigma(dy).
\end{aligned}\end{equation}
Then, for any tempered distribution $U\in\Sa'(\R^d)^d$, we define its normalized flow through the sphere $S(y,r)$ denoted by $\vec\Ra[U]\in\Sa'(\Sigma)$ as

\begin{equation}\nn\begin{aligned}
\left<\vec\Ra[u],\ph\right>_{\Sa'(\Sigma),\Sa(\Sigma)}=\left<u,\vec\Ra^*[\ph]\right>_{\Sa'(\R^d),\Sa(\R^d)},\ \ \ \forall \; \ph\in\Sa(\R^d).
\end{aligned}\end{equation}
\end{de}

The following result is easy to prove. 
\begin{prop} \label{propkl} For any $u\in\Sa'(\R^d)$, $U\in\Sa'(\R^d)^d$, we have the following identities in the sense of distributions:

\begin{equation}\begin{aligned}
\vec \Ra[\g u] = \dr_r\Ra [u],
\end{aligned}\end{equation}

\begin{equation}\begin{aligned}
\vec \Ra[\g\times  U] = 0,
\end{aligned}\end{equation}

\begin{equation}\begin{aligned}
\Ra[\g\p U] = \frac 1 r\dr_r\left(r\vec\Ra[U]\right),
\end{aligned}\end{equation}
and
\begin{equation}\begin{aligned}
\Ra[\La u] = \frac 1 r\dr_r\big(r\dr_r\Ra[u]\big).
\end{aligned}\end{equation}
\end{prop}

\bigskip
\bigskip

\section{Recovering the internal data}

The aim of this section is to recover the internal data $\Psi$ with enough stability in order to use it in the next section to recover the absorption map $a$. The section is divided into five steps.

In the first step, we prove that when $a$ belongs to $SBV^\ii(\Om)$, the approximation

\begin{equation}\nn\begin{aligned}
M_\eta(y,r)=-\frac 1{\eta^2}\io\Ph^2(x)v_{y,r,\eta}(x)\p Da(dx)+\cO\left(\eta^{\frac{d-1}{2d}}\right)
\end{aligned}\end{equation}
holds as $\eta$ goes to zero. In the second step, we link the approximated measurement to $\vec\Ra[\Ph^2Da]$ through the exact formula:

\begin{equation}\nn\begin{aligned}
\frac 1{\eta^2}\io\Ph^2(x)v_{y,r,\eta}(x)\p Da(dx)=\left(\left[\vec\Ra[\Ph^2Da]\right]*\left[r^{d-2}w_\eta\right]\right)(y,r),
\end{aligned}\end{equation}
where $*$ is the convolution product with respect to the variable $r$ and $w_\eta(r)=\frac 1 \eta w(r/\eta)$. In the third step, we give a \index{weak Helmholtz decomposition} weak Helmholtz decomposition of $$\Ph^2Da=D\Psi+\g \times G,$$ where $\Psi\in H^s(D)$ with $s\in[0,1/2[$ and is of class $\C^\ii$ outside of $\supp(Da)$ and satisfies $\Psi|_Y=0$. In the fourth step, we prove that its spherical means Radon transform $\Ra[\Psi]$ is stably approximated in the space $H^{(d-1)/2}(\Sigma)$ in order to satisfy  the assumptions of the Palamodov theorem. We conclude by proving the stable reconstruction of $\Psi$ in $L^2(D)$, where $D$ is a smooth subdomain of interest containing $\supp(Da)$ and is such that $Y\su\Om \setminus \overline{D}$.
\bigskip

\subsection{Step 1: From physical to ideal measurements}

\begin{de}[Ideal measurements] We call the ideal measurement function associated to the absorption $a\in SBV^\ii(\Om)$ the function defined on $\Sigma$ by
\begin{equation}\begin{aligned}
\ti M_\eta(y,r)=-\frac 1 {\eta^2}\int_\Om\Ph^2(x)v_{y,r,\eta}(x)\cdot D a(dx).
\end{aligned}\end{equation} 
\end{de}
In order to prove that $M_\eta$ is close to $\ti M_\eta$ when $\eta$ goes to zero, we need several definitions.

\begin{de}[Wrap condition] Let $\Om'\Subset\Om$ be a smooth domain. We say that the surface $Y \subset \Om \setminus \overline{\Om'}$ satisfies the wrap condition around $\Om'$ if there exists a constant $C>0$ such that for any $x\in \Om'$, $\Gamma\subset S^{d-1}$ measurable, we have
\begin{equation}\nn\begin{aligned}
\sigma\left(Y\cap \mathrm{Cone}(x,\Gamma)\right)\leq C\ \sigma\left(\Gamma\right),
\end{aligned}\end{equation}
where $\mathrm{Cone}(x,\Gamma)=\{x+t\xi,\ \xi\in\Gamma,\ t\in\R^+\}$.
\end{de}

\begin{theo}\label{theo_meas}
Let $a\in SBV^\ii(\Om)$ and let $\Om'$ be such that $\dist(\Om',Y)\geq r_0>0$. Suppose that $Y$ satisfies the wrap condition around $\Om'$. Then, there exists a constant $C>0$ depending on $\Om$, $\Ph$, $|Y|$, $|D a|(\Om)$, $r_0$ and the wrap constant such that

\begin{equation}\nn\begin{aligned}
\n{M_\eta-\ti M_\eta}{L^2(\Sigma)}\leq C\eta^{\frac {d-1}{2d}},
\end{aligned}\end{equation}
and

\begin{equation}\nn\begin{aligned}
\n{P[M_\eta]-P[\ti M_\eta]}{H^1(\Sigma)}\leq C\eta^{\frac {d-1}{2d}},
\end{aligned}\end{equation}
where $P$ is the operator defined by

\begin{equation}\nn\begin{aligned}
P[\ph](y,r)=- \int_0^r\frac{\ph(y,\rho)}{\rho^{d-2}}d\rho.
\end{aligned}\end{equation}
\end{theo}

To prove this result, we need several lemmas. The first one is a spherical density result for the Radon measure $|Da|$. Its proof uses some measure density results and is given in Appendix \ref{app_density}.

\begin{lem}\label{lem_ph} Consider $a\in SBV^\ii(\Om)$ constant out of the subdomain $D\Subset \Om$ and let the mollifier sequence $w_\eta(r)=\frac 1 \eta w\left(\frac 1\eta\right)$, where $w$ is given by (\ref{eq_defv}). Suppose that $Y$ satisfies the wrap condition around $D$. Then, the sequence of functions defined on $\Sigma$ by

\begin{equation}\nn\begin{aligned}
\ph_\eta(y,r)=\int_\Om w_\eta(|x-y|-r)|D a|(dx)
\end{aligned}\end{equation}
satisfies

\begin{equation}\nn\begin{aligned}
\n{\ph_\eta}{L^2(\Sigma)}\leq C\eta^{-\frac 1{2d}}
\end{aligned}\end{equation}
with $C$ depending on $|D a|(\Om)$, $|Y|$, and  the wrap constant. 
\end{lem}

\medskip
In the next lemma, we rewrite the measurement map $M_\eta$.

\begin{lem}\label{lem_rewrite} For any $(y,r)\in\Sigma$, we have
\begin{equation}\nn\begin{aligned}
M_\eta(y,r)=-\frac 1{\eta^2}\int_\Om T[v_{y,r,\eta}](x,y)v_{y,r,\eta}(x)\cdot  Da(dx),  
\end{aligned}\end{equation}
where
\begin{equation}\nn\label{eq_f}\begin{aligned}
T[v](x,y)=\int_0^1(\Ph\Ph_v)\left(x+tv(x)\right)\left(1+t\frac{|v(x)|}{|x-y|}\right)^{d-1}dt.
\end{aligned}\end{equation}
\end{lem}

\Proof Since we fix $y$ (supposed to be zero), $r>r_0$ and $\eta>0$, we will not write the dependence with respect to these variables. We first introduce an approximation sequence of smooth functions $(a^\e)_{\e>0}$ such that $\supp(a^\e-a_0)\subset\Om'$ and $a^\e \rightarrow a \tin L^2(\Om)$. Note that its derivative $ \g a^\e$ converges to $ Da$ for the $H^{-1}(\Om)^d$ norm.

We define now a flow $\ph(x,t)=x+tv(x)$, $\ph\in\C^\ii\big(\R^d\times[0,1],\R^d\big)$. The condition $w'>-1$ ensures that this flow is invertible in the sense that there exists a flow $\ph^{-1}(x,t)$ of class $\C^\ii$ such that $\ph(\ph^{-1}(x,t),t)=\ph^{-1}(\ph(x,t),t)=x$ for all $(x,t)\in\R^d\times [0,1]$. In particular, it satisfies for any $x\in\R^d$, $\ph^{-1}(x,0)=x$ and $\ph^{-1}(x,1)=(Id+v)^{-1}(x)$. For all $x\in\R^d$, $\e>0$, we have
\begin{equation}\nn\begin{aligned}
a^\e\circ(Id+v)^{-1}(x)-a^\e(x) &=\int_0^1 \g a^\e(\ph^{-1}(x,t))\cdot \dr_t\ph^{-1}(x,t)dt\\
\int_\Om(a^\e_v-a^\e)p &=\int_\Om\int_0^1 \g a^\e(\ph^{-1}(x,t))\cdot \dr_t\ph^{-1}(x,t)p(x)dtdx\\
\int_\Om(a^\e_v -a^\e)p &=\int_0^1\int_\Om \g a^\e(\ph^{-1}(x,t))\cdot \dr_t\ph^{-1}(x,t)p(x)dxdt,
\end{aligned}\end{equation}
where $p= \Phi \Phi_v$. 

Hence, using the change of variables $x\mapsto\ph(x,t)$, we get
\begin{equation}\nn\begin{aligned}
\int_\Om(a^\e_{v}-a^\e)p &=\int_0^1\int_\Om \g a^\e(x)\cdot \dr_t\ph^{-1}(\ph(x,t),t)p\circ\ph(x,t)\det(d_x\ph(x,t))dxdt\\
&=-\int_\Om F\cdot \g a^\e,
\end{aligned}\end{equation}
where

\begin{equation}\nn\begin{aligned}
F(x)=-\int_0^1\dr_t\ph^{-1}(\ph(x,t),t)p\circ\ph(x,t)\det(d_x\ph(x,t))dt.
\end{aligned}\end{equation}
As $p\in H^2(\Om)$, the function $F$ belongs to $H^1(\Om)^d$. Passing to the limit when $\e$ goes to zero in the previous equation, the term in the left-hand side goes to $\int_\Om(a_u-a)p$ and as $F\in H^1(\Om)^d$ and $\supp\big( \g a^\e\big)\subset\Om'\subset\subset\Om$, the right-hand side converges to $\int_\Om F(x)\cdot  Da(dx)$. Hence, the formula
\begin{equation}\nn\begin{aligned}
M=\int_\Om(a_v-a)p=-\int_\Om F(x)\cdot  Da(dx)
\end{aligned}\end{equation}
holds. In order to simplify the writing of $F$, we recall two useful properties satisfied by $\ph$ and $\ph^{-1}$. Deriving the identity $\ph^{-1}(\ph(x,t),t)=x$ with respect to $t$ and $x$, we get
\begin{equation}\nn\begin{aligned}
d_x\ph^{-1}(\ph(x,t),t)\dr_t\ph(x,t)+\dr_t\ph^{-1}(\ph(x,t),t) &=0,\\
d_x\ph^{-1}(\ph(x,t),t)d_x\ph(x,t) &=Id.
\end{aligned}\end{equation}
We recall that $d_x\ph(x,t)=Id+tdv(x)$. Now noticing that $\dr_t\ph^{-1}(\ph(x,t),t)=-(Id+tdv(x))$, we rewrite $F$ as follows:
\begin{equation}\nn\begin{aligned}
F(x)=-\int_0^1p\big(x+tv(x)\big)\det\big(Id+tdv(x)\big)\big(Id+tdv(x)\big)^{-1}v(x)dt.
\end{aligned}\end{equation}
Fortunately, $dv(x)$ is diagonal in the spherical orthonormal basis $\B=(\xi,e_2,\cdots,e_d)$, where $\xi={x}/{|x|}$ and $(e_2,\cdots,e_d)$ is an orthonormal basis of $\xi^\perp$, the hyperplane orthogonal to $\xi$. Its matrix in this basis is given by
\begin{equation}\nn\begin{aligned}
\mat_\B(dv(x))=\left[\begin{matrix} \frac{r_0}{r}w'\left(\frac{|x|-r}{\eta}\right) & 0\\  0 & \frac{|v(x)|}{|x|}I_{d-1}\end{matrix}\right].
\end{aligned}\end{equation}
Then,
\begin{equation}\nn\begin{aligned}
\mat_\B(Id+tdv(x)) &=\left[\begin{matrix} 1+t\frac{r_0}{r}w'\left(\frac{|x|-r}{\eta}\right) & 0\\  0 & \left(1+t\frac{|v(x)|}{|x|}\right)I_{d-1}\end{matrix}\right],
\end{aligned}\end{equation}
and from this matrix we deduce that
\begin{equation}\nn\begin{aligned}
\det\big(Id+tdv(x)\big) &=\left[1+t\frac{r_0}{r}w'\left(\frac{|x|-r}{\eta}\right)\right]\left[1+t\frac{|v(x)|}{|x|}\right]^{d-1}\\
\big(Id+tdv(x)\big)^{-1}v(x) &=\frac{v(x)}{1+t\frac{r_0}{r}w'\left(\frac{|x|-r}{\eta}\right)}.
\end{aligned}\end{equation}
Therefore, 
\begin{equation}\nn\begin{aligned}
F(x) &=\int_0^1p\big(x+tv(x)\big)\left[1+t\frac{|v(x)|}{|x|}\right]^{d-1}dt\ v(x).
\end{aligned}\end{equation}
Replacing $|x|$ by $|x-y|$ and rewriting the dependence in $y,r$, and $\eta$, we finally get the expected formula.\cqfd

The next result shows that the shifted absorption map $a_v$ stays close to $a$ in $L^1(\Om)$ if $\eta$ is small. The key result is optimal in the sense that it requires that $a$ to be of bounded variation. In fact, it shows that any reconstruction would be impossible without this minimal regularity.

\begin{prop}\label{prop_av-a} Consider $a\in \A_0\cap BV(\Om)$ and let the internal displacement $v$ be given by (\ref{eq_defv}). We have the following estimate:

\begin{equation}\nn\begin{aligned}
\n{a_{v}-a}{L^1(\Om)}\leq C| Da|(\Om)\eta
\end{aligned}\end{equation}
with $C$ depending only on the space dimension $d$.
\end{prop}
\Proof Let us consider an approximation sequence $(a^\e)_{\e>0}\su\C^0(\Om)$ such that $\supp(a^e-a)\su D$ and $\n{a^\e-a}{L^1(\Om)}\leq\e$. Now, we define the flow $\ph\in\C^\ii\big(\R^d\times[0,1]\big)$ by $\ph(x,t)=x+tv_\eta(x)$. The condition $w^\prime > -1$ ensures that this flow is invertible in the sense that there exists a flow $\ph^{-1}(x,t)$ of class $\C^\ii$ such that $\ph(\ph^{-1}(x,t),t)=\ph^{-1}(\ph(x,t),t)=x$ for all $(x,t)\in\R^d\times [0,1]$. In particular, it satisfies for any $x\in\R^d$, $\ph^{-1}(x,0)=x$ and $\ph^{-1}(x,1)=(Id+v_\eta)^{-1}(x)$.

For all $x\in\R^d$, $\e>0$, we get

\begin{equation}\nn\begin{aligned}
a^\e_{v_\eta}(x)-a^\e(x) &=a^\e\circ\ph^{-1}(x,1)-a^\e\circ\ph^{-1}(x,0)=\int_0^1 \g a^\e\circ\ph^{-1}(x,t)\dr_t\ph^{-1}(x,t)dt\\
\n{a^\e_{v_\eta}(x)-a^\e(x)}{L^1(\Om)}&\leq \int_0^1\int_\Om\left| \g a^\e\circ\ph^{-1}(x,t)\cdot\dr_t\ph^{-1}(x,t)\right|dxdt\\
&\leq \int_0^1\int_\Om\left| \g a^\e(x)\cdot\dr_t\ph^{-1}(\ph(x,t),t)\right||\det d_x\ph(x,t)|dxdt.
\end{aligned}\end{equation}
A similar computation to the one in the proof of (\ref{lem_rewrite}) leads to $$|\dr_t\ph^{-1}(\ph(x,t),t)\det d_x\ph(x,t)|\leq \left(1+\frac{d(d-1)}{2}\right)|v_\eta(x)|$$ and so,

\begin{equation}\nn\begin{aligned}
\n{a^\e_{v_\eta}(x)-a^\e(x)}{L^1(\Om)}&\leq \left(1+\frac{d(d-1)}{2}\right)\eta\int_\Om| \g a^\e|.
\end{aligned}\end{equation}
Passing now to the limit when $\e$ goes to zero, we get the expected result.\cqfd

As a consequence of Proposition \ref{prop_av-a}, we deduce that the modified light fluence $\Ph_v$ is close to $\Ph$ in $\htwo$ when $\eta$ is small.

By combining (\ref{prop_av-a}) and (\ref{lem_Fcont}), the following result holds.
\begin{cor} Consider $a\in\A_0 \cap BV(\Om)$ and the internal displacement $v$ given by (\ref{eq_defv}). We have the following estimate:

\begin{equation}\nn\begin{aligned}
\n{\Ph_{v}-\Ph}{H^2(\Om)}\leq C\ol \Ph(\ol a-\ul a)^\frac{1}{2}| Da|(\Om)^\frac{1}{2}\eta^\frac{1}{2},
\end{aligned}\end{equation}
where $C$ depends on $d$ and $\Om$.
\end{cor}

\begin{lem}\label{lem_fphi} Consider a subdomain $\Om'\su\Om$ such that $\dist(\Om',Y)\geq r_0>0$. There exists a constant $C>0$ depending on $\Om$, $\Ph$ and $a$ such that
\begin{equation}\nn\begin{aligned}
\n{T[v_{y,r,\eta}](.,y)-\Ph^2}{L^\ii(\Om')}\leq C\eta^{\frac 1 2}.
\end{aligned}\end{equation}
\end{lem}

\Proof For fixed $\eta>0$ and $(y,r)\in \Sigma$, for $t\in[0,1]$ and $x\in\Om'$,

\begin{equation}\nn\begin{aligned}
|(\Ph\Ph_v)(x+tv(x))-\Ph^2(x)| &\leq |\Ph^2(x+tv(x))-\Ph^2(x)|+|(\Ph\Ph_v)(x+tv(x))-\Ph^2(x+tv(x))|\\
&\leq 2\ol\Ph|\Ph(x+tv(x))-\Ph(x)|+\ol\Ph|(\Ph_v(x+tv(x))-\Ph(x+tv(x))|\\
&\leq 2\ol\Ph\n\Ph{\C^{0,\frac 1 2}\left(\ol \Om\right)} \eta^{\frac 1 2}+\ol\Ph C_1\eta^{\frac 1 2}\\
&\leq C\eta^{1/2}.\\
\end{aligned}\end{equation}
Recalling that for $x\in\Om'$, $|x-y|\geq r_0$, we use the previous inequality in (\ref{eq_f}) to get the desired result.\cqfd

Now, we are ready to prove Theorem \ref{theo_meas}.

\Proof (of Theorem \ref{theo_meas}) For any $(y,r)\in \Sigma$,

\begin{equation}\nn\begin{aligned}
|M_\eta-\ti M_\eta|(y,r)&\leq\int_\Om|T[v_{y,r,\eta}](x,y)-\Ph^2(x)|\frac{|v_{y,r,\eta}(x)|}{\eta^2}| Da|(dx)\\
&\leq \n{T[v_{y,r,\eta}]-\Ph^2}{L^\ii(\Om')}\int_\Om w_\eta(|x-y|-r)| Da|(dx).
\end{aligned}\end{equation}
Applying Lemmas \ref{lem_fphi} and \ref{lem_ph}, we get the first inequality,

\begin{equation}\nn\begin{aligned}
\n{M_\eta-\ti M_\eta}{L^2(\Sigma)}\leq C\eta^{\frac {d-1}{2d}}.
\end{aligned}\end{equation}
Next, taking the derivative with respect to the variable $y$,
it follows that
\begin{equation}\nn\begin{aligned}
d_y(M_\eta-\ti M_\eta)(y,r)=\int_\Om d_y \left( \left(T[v_{y,r,\eta}](x,y)-\Ph^2(x)\right)\frac{v_{y,r,\eta}(x)}{\eta^2}\right)\cdot Da(dx)
\end{aligned}\end{equation}
with
\begin{equation}\nn\begin{aligned}
d_y \left( \left(T[v_{y,r,\eta}](x,y)-\Ph^2(x)\right)\frac{v_\eta(x,y,r)}{\eta^2}\right) &=d_vT[v_{y,r,\eta}](x,y)\cdot d_yv_{y,r,\eta}(x) \frac{v_{y,r,\eta}(x)}{\eta^2}\\
&+ \left(T[v_{y,r,\eta}](x,y)-\Ph^2(x)\right)\frac{ d_yv_{y,r,\eta}(x)}{\eta^2}\\
&+ d_yT[v_{y,r,\eta}](x,y) \frac{v_{y,r,\eta}(x)}{\eta^2}.
\end{aligned}\end{equation}
The vector field $v_{y,r,\eta}$ satisfies
\begin{equation}\nn\begin{aligned}
d_yv_{y,r,\eta}=\dr_rv_{y,r,\eta}+|v_{y,r,\eta}|B(y,r,x),
\end{aligned}\end{equation}
where $B(y,r,x)$ is a matrix uniformly bounded with respect to all variables. Moreover,
\begin{equation}\nn\begin{aligned}
d_yT[v_{y,r,\eta}](x,y)=\cO(\eta)
\end{aligned}\end{equation}
with reminder uniform with respect to all variables. Thus we write
\begin{equation}\nn\begin{aligned}
&d_y \left( \left(T[v_{y,r,\eta}](x,y)-\Ph^2(x)\right)\frac{v_\eta(x,y,r)}{\eta^2}\right)\\
 &=\dr_r \left( \left(T[v_{y,r,\eta}](x,y)-\Ph^2(x)\right)\frac{v_{y,r,\eta}(x)}{\eta^2}\right)\frac{(x-y)^T}{|x-y|}\\
&+ R(y,r,x)\frac{|v_{y,r,\eta}(x)|}{\eta^2},
\end{aligned}\end{equation}
where the reminder $R(y,r,x)=\cO(\eta^{1/2})$ uniformly with respect to all variables. 
Here, $T$ denotes the transpose.

Using this identity, we can integrate by parts with respect to $r$ to get
\begin{equation}\nn\begin{aligned}
\int_0^r &\frac 1 {\rho^{d-2}}d_y \left( \left(T[v_\eta(y,\rho)](x,y)-\Ph^2(x)\right)\frac{v_\eta(y,\rho)(x)}{\eta^2}\right)d\rho\\
&=\frac 1 {r^{d-2}}\left( \left(T[v_{y,r,\eta}](x,y)-\Ph^2(x)\right)\frac{v_{y,r,\eta}\g(x)}{\eta^2}\right)\frac{(x-y)^T}{|x-y|}\\
&+(d-2)\int_0^r\frac 1 {\rho^{d-1}}\left( \left(T[v_\eta(y,\rho)](x,y)-\Ph^2(x)\right)\frac{v_\eta(y,\rho)(x)}{\eta^2}\right)d\rho\frac{(x-y)^T}{|x-y|}\\
&+ \int_0^r R(y,\rho,x)\frac{|v_\eta(y,\rho)(x)|}{\eta^2}d\rho.
\end{aligned}\end{equation}
Finally, we integrate over $\Om$ and use Lemma 
\ref{lem_ph}  in order to control these three terms. 
We control all of these by $\eta^{\frac{d-1}{2d}}$.\cqfd
\bigskip

\subsection{Step 2: Linking the measurement with $\Ph^2 Da$}

In the previous subsection, we have shown that the measurement $M_\eta$ is approximated by 

\begin{equation}\nn\begin{aligned}
\ti M_\eta(y,r)=-\frac 1 {\eta^2}\io\Ph^2(x)v_{y,r,\eta}(x)\p Da(dx)
\end{aligned}\end{equation}
in a certain sense when $\eta$ goes to zero. We suggest here another form of  $\ti M_\eta$ using the spherical operators defined in subsection \ref{sdfgd}. As $Da$ is a finite measure compactly supported, it is a tempered distribution on $\R^d$. Since  $\Ph\in \C^0\big(\ol\Om\big)$, the vector field $\Ph^2Da$ is a tempered distribution on $\R^d$ defined by

\begin{equation}\nn\begin{aligned}
\left<\Ph^2Da,\ph\right>_{S'(\R^d)^d,S(\R^d)^d}=\io\Ph^2\ph\cdot Da.
\end{aligned}\end{equation}

The following result holds.
\begin{prop} For any $a\in\A$, $\eta>0$ we have the formula

\begin{equation}\label{eq_M2R}\begin{aligned}
\ti M_\eta=-\frac 1 r\left[\left({r^{d-1}}\vec\Ra[\Ph^2Da]\right)*w_\eta(-.)\right]\quad \tin\Sigma,
\end{aligned}\end{equation}
where $*$ is the one dimensional convolution product with respect to the variable $r$ and $w_\eta(r)=\dfrac 1\eta w\left(\dfrac r\eta\right)$.
\end{prop}

\Proof Consider a test function $\ph\in\Sa(\Sigma)$. We have

\begin{equation}\nn\begin{aligned}
-\int_\Sigma \ti M_\eta \ph&=-\int_Y\int_0^\ii M_\eta(y,r)\ph(y,r)\s(dy)dr\\
&=\frac 1 {\eta^2}\io\Ph^2(x)\left(\int_Y\int_0^\ii v_{y,r,\eta}(x)\ph(y,r)\s(dy)dr\right)\p Da(dx)\\
&=\io\Ph^2(x)\left(\int_Y\int_0^\ii \frac 1 r w_\eta\left(|x-y|-r\right)\ph(y,r)dr\frac{x-y}{|x-y|}\s(dy)\right)\p Da(dx)\\
&=\io\Ph^2(x)\left(\int_Y \left(w_\eta*\frac {\ph(y,.)} r\right)(|x-y|) \frac{x-y}{|x-y|}\s(dy)\right)\p Da(dx)\\
&=\io\Ph^2(x)\vec\Ra^*\left[r^{d-1}\left(w_\eta*\frac {\ph(y,.)} r\right)\right](x)\p Da(dx)\\
&=\left<\Ph^2Da,\vec\Ra^*\left[r^{d-1}\left(w_\eta*\frac {\ph(y,.)} r\right)\right]\right>_{\Sa'(\R^d),\Sa(\R^d)}\\
&=\left<\vec\Ra[\Ph^2Da],r^{d-1}\left(w_\eta*\frac {\ph(y,.)} r\right)\right>_{\Sa'(\Sigma),\Sa(\Sigma)}\\
&=\left<\frac 1 {r^{d-1}}\vec\Ra[\Ph^2Da],w_\eta*\frac {\ph(y,.)} r\right>_{\Sa'(\Sigma),\Sa(\Sigma)}\\
&=\left<\frac 1 r\left({r^{d-1}}\vec\Ra[\Ph^2Da]\right)*w_\eta(-.),\ph\right>_{\Sa'(\Sigma),\Sa(\Sigma)}.\\
\end{aligned}\end{equation}\cqfd
\bigskip

\subsection{Step 3: Helmholtz decomposition of $\Ph^2 Da$}

Since $\Ph^2 Da$ is a tempered distribution, we can consider its Fourier transform. As $a\in H^s(\Om)$ with $s\in[0,1/2[$ and $Da$ is supported in some compact subset $K$ of $\Omega$ , it follows that $Da\in H^{s-1}_K(\Om)^d$ (see Appendix \ref{appendsobolev} for the definition of $H_K^{s - 1}(\Omega)$). Moreover, as $\Ph^2$ is in $\htwo$, we also have that $\Ph^2Da\in H^{s-1}_K(\Om)^d$. From Appendix \ref{appendsobolev}, we deduce that $\wh{\Ph^2Da}$ belongs to $L^1_{\mathrm{loc}}(\R^d)^d$ and satisfies

\begin{equation}\nn\begin{aligned}
\int_{\R^d}\left|\wh{\Ph^2Da}\right|^2(\xi)\left(1+|\xi|^2\right)^{(s-1)}d\xi< +\ii .
\end{aligned}\end{equation}

Let the Sobolev space $H^{\al+1}_{\mathrm{curl}}(\R^d)$ be defined by
$$
H^{\al+1}_{\mathrm{curl}}(\R^d):=\left\{A\in H^\al(\R^d)^d,\ \g\times A\in H^\al(\R^d)^d\right\}.
$$
The following proposition gives a generalization of the Helmholtz decomposition for some compactly supported distributional vector fields.

\begin{prop} Consider $\al\in\R$ and $U\in H^\al_K(\Om)^d$, where $K$ is a compact of $\Om$. There exists $u\in H^{\al+1}(\R^d)$ and $A\in H^{\al+1}_{\mathrm{curl}}(\R^d)$ such that

\begin{equation}\nn\begin{aligned}
U=D u+\g\times A
\end{aligned}\end{equation}
in the sense of distributions.
\end{prop}

\Proof As $U\in H^\al_K(\Om)$, $\wh U\in L_{\mathrm{loc}}^1(\R^d)^d$. We define now $\wh u=\frac{\wh U\cdot\xi}{i|\xi|^2}\in L_{\mathrm{loc}}^1(\R^d)$ and $\wh A=\frac{\wh U\wedge\xi}{i|\xi|^2}\in L_{\mathrm{loc}}^1(\R^d)^d$. We have the decomposition $\wh U=i\wh u\xi+i\xi\wedge\wh A$. As $i\wh u\xi$ is the Fourier transform of $\g u$ where $u$ is the inverse Fourier transform of $\wh u$ and has the same integrability as $\wh U$, we deduce that $u\in H^{\al+1}(\R^d)$ and $A$, the inverse Fourier transform of $\wh A$, is in $H^{\al+1}_{\mathrm{curl}}(\R^d)$.\cqfd

Using this last result, we write

\begin{equation}\begin{aligned}
\Ph^2Da=D\Psi+\g\times G
\end{aligned}\end{equation}
with $\psi\in H^s(\R^d)$ and therefore, we get 

\begin{equation}\begin{aligned}
\vec\Ra[\Ph^2\g a]=\dr_r\Ra[\psi]
\end{aligned}\end{equation} in the sense of distributions. Using this last identity in (\ref{eq_M2R}) we obtain that

\begin{equation}\label{eq_Mti}\begin{aligned}
\ti M_\eta=-\frac 1 r\left[\left({r^{d-1}}\dr_r\Ra[\Psi]\right)*w_\eta \right]\quad \tin\Sigma.
\end{aligned}\end{equation}
Equation \eqref{eq_Mti} plays an important role to recover the internal data in the next step.
\bigskip

\subsection{Step 4: Approximating $\Ra[\Psi]$}

We can show now that from the previous identity, the quantity $\Ra[\Psi]$ can be approximated up to a function depending only in $y$ in $H^{\frac{d-1}2}(\Sigma)$ in order to apply Palamodov's theorem (see \cite{victoire}).

\begin{theo}\label{theo_step4} Let $s\in\left]\frac 1 3, \frac 1 2\right[$ and $\al>0$ and consider $a\in\A_{s+\al}$. Then $P[\ti M_\eta]$ converges to $\Ra[\Psi]-g$ in $H^{\frac{d-1}2+s}(\Sigma)$ where $g$ is a function depending only on $y$. More precisely, there exists a constant $C$ depending on $d$, $s$, $\al$, and  $\Sigma$ such that
\begin{equation}\nn\begin{aligned}
\n{P[\ti M_\eta]-\Ra[\Psi]+g}{H^{\frac{d-1}2+s}(\Sigma)}\leq C\eta^\frac{\al}{\al+1}\n{\Psi}{H^{s+\al}(\Sigma)}.
\end{aligned}\end{equation}
\end{theo}

\Proof Starting from (\ref{eq_Mti}) and integrating by parts, we write

\begin{equation}\nn\begin{aligned}
\ti M_\eta(y,r)=\int_\R\Ra[\Psi](y,\rho)\dr_\rho\left(\frac{1}{r}w_\eta(\rho-r)\rho^{d-1}\right)d\rho.
\end{aligned}\end{equation}
Now, applying $P$ to $\ti M_\eta$, we get

\begin{equation}\nn\begin{aligned}
P[\ti M_\eta](y,r)=-\int_\R\Ra[\Psi](y,\rho)\dr_\rho\left(\int_{r_0}^rw_\eta(\rho-s)\frac{\rho^{d-1}}{s^{d-1}}ds\right)d\rho.
\end{aligned}\end{equation}
Let us develop the test function as follows:

\begin{equation}\nn\begin{aligned}
\dr_\rho\int_{r_0}^rw_\eta(\rho-s)\frac{\rho^{d-1}}{s^{d-1}}ds=w_\eta(\rho-r_0)\frac{\rho^{d-1}}{r_0^{d-1}}-w_\eta(\rho-r)\frac{\rho^{d-1}}{r^{d-1}}-\theta_\eta(\rho,r)
\end{aligned}\end{equation}
with
\begin{equation}\nn\begin{aligned}
\theta_\eta(\rho,r)=(d-1)\int_{r_0}^rw(\rho-s)(s-\rho)\frac{\rho^{d-2}}{s^{d}}ds,
\end{aligned}\end{equation}
which satisfies $\n{\theta_\eta}{H^1\left(]0,R[^2\right)}\leq C\eta^{1/2}$, where $C$ depends on $r_0$, $R$, and $d$. Finally, we write

\begin{equation}\nn\begin{aligned}
P[\ti M_\eta](y,r) &=\frac 1{r^{d-1}}\int_0^R\rho^{d-1}\Ra[\Psi](y,\rho)w_\eta(\rho-r)d\rho-\frac 1{r_0^{d-1}}\int_0^R\rho^{d-1}\Ra[\Psi](y,\rho)w_\eta(\rho-r_0)d\rho\\
 &+\int_0^R\Ra[\Psi](y,\rho)\theta_\eta(\rho,r)d\rho.
\end{aligned}\end{equation}
Using Lemma \ref{lem_ker1}, we bound the $H^{\frac{d-1}2}(\Sigma)$ norm of the third term by $C\eta^{\frac 1 2}\n{\Ra[\Psi]}{H^{\frac{d-1}2}(\Sigma)}$. Moreover, using Lemma \ref{lem_ker2}, we say that as $\Ra[\Psi]\in H^{\frac{d-1}2+s+\al}(\Sigma)$, the first term converges to $\Ra[\Psi]$ for the norm of $H^{\frac{d-1}2+s}(\Sigma)$ with an error controlled by $\eta^\frac{\al}{\al+1}\n{\Ra[\Psi]}{H^{\frac{d-1}2+s+\al}(\Sigma)}$. Using the same argument, the second term goes to $g(y)=\Ra[\Psi](y,r_0)$ in the same manner. We finally obtain that

\begin{equation}\nn\begin{aligned}
\n{P[\ti M_\eta]-\Ra[\Psi]+g}{H^{\frac{d-1}2+s}(\Sigma)}\leq C\eta^\frac{\alpha}{\alpha+1}\n{\Ra[\Psi]}{H^{\frac{d-1}2+s+\al}(\Sigma)},
\end{aligned}\end{equation}
where $C$ depends on $d$ and the manifold $\Sigma$. \cqfd
\bigskip

\subsection{Step 5: Approximating $\Psi$}

We recall here that we have assumed the invertibility of the spherical means Radon transform. We apply $\Ra^{-1}$ to the inequality given in Theorem \ref{theo_step4} to get

\begin{equation}\begin{aligned}
\n{\Ra^{-1}\circ P[\ti M_\eta]-\Psi+\Ra^{-1}[g]}{H^s(D)}\leq C \eta^\frac{\al}{\al+1}\n{\Psi}{H^{s+\al}(\Sigma)}
\end{aligned}\end{equation}
for some positive constant $C$ independent of $\eta$. 

The problem that we have here is that we do not know the map $h=\Ra^{-1}[g]$. 
Nevertheless, if we assume that $Y = \partial D'$
with $D \Subset D'$ being smooth, then, 
since $Da$ vanishes outside of $D$ 
and $Y \subset \Om \setminus \overline{D}$, we can find $\ti\Psi \in H^s$ such that 
$\ti\Psi|_Y=0$ and the Helmholtz decomposition \cite{helmdec}
$$
\Ph^2Da =D\ti\Psi+\g\times\ti G
$$
holds in $D'$. Extending $\ti\Psi$ by zero outside $D'$ yields an $H^s$ function for $s \in ]0,1/2[$. 
Because $\dr_r\Ra[h]=0$ and using Proposition \ref{propkl}, we have $\Ra[\La h]=0$ and therefore $\La h=0$ in $D'$. The boundary condition $\ti\Psi|_Y=0$ shows that $h$ is uniquely determined by solving a Dirichlet problem for the Laplacian. 

%
%
%

%

To summarize five steps of this section, we state the following theorem.

\begin{theo}\label{theo_recapsection2} Consider $a\in\A_0\cap SBV^\ii(\Om)$ satisfying $\supp(Da)\su D\Subset \Om$ and assume that $Y$ satisfies the wrap condition around $D$ and surrounds $D$. Then there exists $\Psi\in H^s(\R^d)$ with $s\in]0,1/2[$ and of class $\C^\ii$ outside of $\supp(Da)$ satisfying $\Psi|_Y=0$ and

\begin{equation}\label{eq_recapsection2}\begin{aligned}
\Ph^2Da=D\Psi + \g\times G,
\end{aligned}\end{equation}
where $G \in H^{s}_{\mathrm{curl}}(\R^d)$ and $\Ph=F[a]$. Moreover, $\Ra^{-1}\circ P[M_\eta] +h$ converges strongly to $\Psi$ in $H^s(D)$ when $\eta$ goes to zero at a speed bounded by $\cO\left(\eta^\frac{1}{4}\right)$. Here, $h$ is determined as above.  
\end{theo}

This map $\Psi$ will be now the starting point of the reconstruction procedure. In the next section, we assume that $\Psi$ is known in $H^s(D)$ up to a small error in $H^s(D)$. We will see how to approximate the absorption parameter $a$ from this data.
\bigskip
\bigskip

\section{Stable reconstruction of the absorption map}

In this section, we assume that the assumptions of Theorem \ref{theo_recapsection2} are satisfied and suppose in addition that $Y=\dr D$. This is possible since $Da$ is assumed to be compactly supported in $D$. It simply suffices to enlarge $D$. 
As a consequence, we assume the knowledge of $\Psi\in H^s(D)$ of class $\C^\ii$ in a neighborhood of $\dr D$, which satisfies $\Psi|_{\dr D}=0$. The goal of this section is to present a method to estimate the absorption map $a$ from the knowledge of $\Psi$. We choose $s$ such that $H^{s+1}(\Om)$ is embedded in $L^\ii(\Om)$ in dimensions $2$ and $3$. This is true for any $s\in]1/3,1/2[$. 

Let us take the divergence of \eqref{eq_recapsection2} in the sense of distributions to get

\begin{equation}\nn\begin{aligned}
\g\p(\Ph^2Da)=\La\Psi,
\end{aligned}\end{equation}
which looks like an elliptic equation with unknown $a$. There are two difficulties here. The first one is that we do not have enough regularity to deal with this equation using a variational approach. To do so, we should have $\Psi$ in $H^1(D)$ and look for a solution $a$ in $H^1(D)$. The second difficulty is that the diffusion term $\Ph^2$ is unknown here and depends on $a$ by $\Ph^2= F[a]^2$, where $F$ is the light fluence operator.

Finally, we recall the definition of the set of admissible absorption distributions:
$$
\A_s=\left\{a\in\A_0\cap H^s(\Om),\ \n{a}{H^s(\Om)}\leq R_{\A_s}\right\}
$$
and define 
$$
\B_s = \left\{\Phi \in W^{1,\infty}(D),\ \underline{\Phi} \leq \Phi \leq \ol{\Phi}, 
\| \nabla \Phi\|_{L^\infty} \leq R_{\B_s}
 \right\},
 $$
where $R_{\B_s}= \sup_{a \in \A_s} \| \nabla F[a]\|_{L^\infty(D)}$. Note that $F$ maps $\A_s$ into $\B_s$.  

\bigskip

\subsection{The change of function argument}
The main idea is to introduce a new variable:

\begin{equation}\label{eq_ati}\begin{aligned}
\ti a=a-a_0-\frac{\Psi}{\Ph^2},
\end{aligned}\end{equation}
which is well defined in $H^s(D)$ since $\Ph\geq\underline \Ph$.

\begin{prop} For all $a \in \A_s$ and $\Phi=F[a]$, we have $\ti a\in H^1_0(D)$ . \end{prop}

\Proof In the sense of distributions, we have

\begin{equation}\nn\begin{aligned}
D\ti a &=D a-\frac{D\Psi}{\Ph^2}+2\Psi\frac{\g\Ph}{\Ph^3},\\
\Ph^2D\ti a &=\Ph^2D a-{D\Psi}+2\Psi{\g\log\Ph},\\
\g\p(\Ph^2D\ti a) &=\g\p(2\Psi{\g\log\Ph}),\\
\Ph^2\La\ti a &=\g\p(2\Psi{\g\log\Ph})-\g(\Ph^2)\p D\ti a,\\
\La\ti a &=\frac 1{\Ph^2}\g\p(2\Psi{\g\log\Ph})-2\g(\log\Ph)\p D\ti a.
\end{aligned}\end{equation}
Consider a test function $\ph\in\D(D)$ and using the fact that $\g\Phi\in L^\ii(D)$, which follows from the fact that $\Ph\in H^{2+s}(\Om)$, we have

\begin{equation}\nn\begin{aligned}
&\left<\frac 1{\Ph^2}\g\p(2\Psi{\g\log\Ph}),\ph\right>_{\D'(D),\D(D)} =\left<\g\p(2\Psi{\g\log\Ph}),\frac\ph{\Ph^2}\right>_{H^{-1}(D),H^1_0(D)}\\
&=-2\int_D\frac{\Psi}{\Ph^2}\g(\log\Ph)\p\left({\g\ph}-2\ph{\g(\log\Ph)}\right)\\
&\leq \frac2{\ul\Ph^3}\n{\Psi}{L^2(D)}\n{\g\Phi}{L^\ii(D)}\left(\n{\g\ph}{L^2(D)}+\frac2{\ul\Ph}\n{\g\Phi}{L^\ii(D)}\n{\ph}{L^2(D)}\right)\\
&\leq C\n{\ph}{H^1_0(D)}
\end{aligned}\end{equation}
and so $\frac 1{\Ph^2}\g\p(2\Psi{\g\log\Ph})\in H^{-1}(D)$. We also have

\begin{equation}\nn\begin{aligned}
&\left<2\g(\log\Ph)\p D\ti a,\ph\right>_{\D'(D),\D(D)} =\left<D\ti a,2\g(\log\Ph)\ph\right>_{H^{-1}(D),H^1_0(D)}\\
&=-2\int_D\ti a\big(\ph\La(\log\Ph)+\g(\log\Ph)\p\g\ph\big)\\
&=-2\int_D\ti a\left(a\ph-\frac{|\g\Ph|^2}{\Ph^2}\ph+\g(\log\Ph)\p\g\ph\right)\\
&\leq 2\n{\ti a}{L^2(D)}\left( \ol a\n{\ph}{L^2(D)}+\frac{\n{\g\Phi}{L^\ii(D)}^2}{\ul\Ph^2}\n{\ph}{L^2(D)}+\frac{\n{\g\Phi}{L^\ii(D)}}{\ul\Ph}\n{\g\ph}{L^2(D)}\right)\\
&\leq C\n{\ph}{H^1_0(D)}
\end{aligned}\end{equation}
and so $2\g(\log\Ph)\p D\ti a\in H^{-1}(D)$. Finally, since $\La\ti a\in H^{-1}(D)$ and $\ti a$ is smooth in a neighborhood of $\dr D$ and satisfies $\ti a|_{\dr D}=0$, it follows from the standard regularity theory  that $\ti a\in H^1_0(D)$.\cqfd

From the previous computation, it follows that $\ti a$ is defined as the unique solution of

\begin{equation}\label{eq_G}\left\{\begin{aligned}
\g\p(\Ph^2\g \ti a) &=\g\p(2\Psi{\g\log\Ph})\bb &\tin D,\\
\ti a &=0\bb &\ton\dr D.
\end{aligned}\right.\end{equation}
This system allows us to define an operator

\begin{equation}\begin{aligned}
\ti G_\Psi:\B_s &\ra H^1_0(\Om)\\
\Ph &\mt\ti a
\end{aligned}\end{equation}
and the one which gives $a$ from $\Ph$,

\begin{equation}\begin{aligned}
G_\Psi:\B_s &\ra H^s(\Om)\\
\Ph &\mt \left\{\begin{aligned} &a_0+\ti G_\Psi[\Ph]+\frac{\Psi}{\Ph^2}\ &\tin D, \\&a_0\ &\tin \Om\bs\ol D .\end{aligned}\right.
\end{aligned}\end{equation}
%

The global problem that we have to solve now is to find a pair $(\ti a,\Ph)\in H^1(D)\times\B_s$ such that

\begin{equation}\nn\left\{\begin{aligned}
-\La\Ph+\left(a_0+\ti a+\frac{\Psi}{\Ph^2}\1_D\right)\Ph &=0\ \quad \tin \Om,\\
\g\p(\Ph^2\g \ti a) &=\g\p(2\Psi{\g\log\Ph})\ \quad \tin D,\\
l\dr_\nu\Ph+\Ph &=g\ \quad \ton\dr\Om,\\
\ti a &= 0\ \quad \ton\dr D,\\
\ti a &= 0\ \quad \tin\dr \Om\bs\ol D,\\
\end{aligned}\right.\end{equation}
where $\1_D$ denotes the characteristic function of $D$. 
\bigskip

\subsection{Fixed point algorithm}

%
We look for a solution $a$ as the fixed point of the map $G_\Psi\circ F:\A_s\ra H^s(\Om)$. In order to cycle this operator, we introduce the truncation operator

\begin{equation}\nn\begin{aligned}
T:H^s(\Om) &\ra H^s(\Om)\\
a &\mt \max\big(\min(a,\ol a),\ul a\big)
\end{aligned}\end{equation}
and look for a fixed point of the operator $T\circ G_\Psi\circ F:\A_s\ra H^s(\Om)$.

\begin{theo}\label{theo_contraction} Consider $\Psi$ in $H^s(D)$. The operator $T\circ G_\Psi\circ F:\A_s\ra H^s(\Om)$ is $H^s$-Lipschitz and, for any $a$, $a'\in\A$, we have

\begin{equation}\nn\begin{aligned}
\n{T\circ G_\Psi\circ F}{Lip\left(H^s(\Om)\right)}\leq c(s,d,\Om,\ul\Phi,\ol\Phi,R_{\B_s})\n{\Psi}{H^s(\Om)}
\end{aligned}\end{equation}
and if $\n{\Psi}{H^s(\Om)}$ is small enough, $T\circ G_\Psi\circ F$ is a contraction from $\A_s$ into $\A_s$ and admits a unique fixed point in $\A_s$ called $a_\Psi$.
\end{theo}

\Proof Reconsidering the Lipschitz estimate of the system \ref{eq_light} with $a\in H^s(\Om)$ and taking into account that $F[a]$ belongs to $W^{1,\ii}(\Om)$ (see the last paragraph of Subsection \ref{subsection 2.2}) we can deduce that $F:H^s(\Om)\ra W^{1,\ii}(D)$ is Lipschitz and obtain that

\begin{equation}\nn\begin{aligned}
\n{F}{Lip\left(H^s(\Om),W^{1,\ii}(D)\right)}\leq c(s,d,\Om)\left(R_{\B_s}+\ol\Ph\right).
\end{aligned}\end{equation}
Consider now $\Ph$ and $\Ph'$ in $\B_s$, then

\begin{equation}\nn\begin{aligned}
|\g\log\Ph-\g\log\Ph'|\leq \frac 1 {\ul\Ph}|\g(\Ph-\Ph')|+\frac {|\g\Ph'|} {\ul\Ph^2}|\Ph-\Ph'|,
\end{aligned}\end{equation}
and so

\begin{equation}\nn\begin{aligned}
\n{\g\log\Ph-\g\log\Ph'}{L^\ii(\Om)}\leq \frac 1 {\ul\Ph}\left(1+\frac{R_{\B_s}}{\ul\Ph}\right)\n{\Ph-\Ph'}{W^{1,\ii}(\Om)}.
\end{aligned}\end{equation}
This inequality proves that $\ti G_\Psi:\B_s \ra H^1_0(\Om)$ is Lipschitz and that

\begin{equation}\nn\begin{aligned}
\n{\ti G_\Psi}{Lip\left(W^{1,\ii}(D),H^1_0(\Om)\right)}\leq\frac 1 {\ul\Ph^3}\left(1+\frac{R_{\B_s}}{\ul\Ph}\right)\n{\Psi}{L^2(D)}.
\end{aligned}\end{equation}
We can now control the Lipschitz norm of $G_\Psi$. Noticing that

\begin{equation}\nn\begin{aligned}
\n{\frac 1 {\Ph^2}-\frac 1 {\Ph'^2}}{W^{1,\ii}(D)}\leq\frac 1 {\ul\Ph^3}\left(2+3\frac{R_{\B_s}}{\ul\Ph}\left(\frac{\ol\Ph}{\ul\Ph}\right)^2 \right)\n{\Ph-\Ph'}{W^{1,\ii}(D)},
\end{aligned}\end{equation}
we get 

\begin{equation}\nn\begin{aligned}
\n{G_\Psi\Ph-G_\Psi\Ph'}{H^s(\Om)} &\leq \n{\ti G_\Psi\Ph-\ti G_\Psi\Ph'}{H^1_0(\Om)}+\n{\Psi}{H^s(\Om)}\n{\frac 1 {\Ph^2}-\frac 1 {\Ph'^2}}{W^{1,\ii}(D)}\\
&\leq \frac 1 {\ul\Ph^3}\left[3+3\frac{R_{\B_s}}{\ul\Ph}\left(1+\left(\frac{\ol\Ph}{\ul\Ph}\right)^2\right) \right]\n{\Psi}{H^s(\Om)}\n{\Ph-\Ph'}{W^{1,\ii}(D)}
\end{aligned}\end{equation}
and finally,

\begin{equation}\nn\begin{aligned}
\n{G_\Psi}{Lip\left(W^{1,\ii}(D),H^s(\Om)\right)}\leq \frac 1 {\ul\Ph^3}\left[3+3\frac{R_{\B_s}}{\ul\Ph}\left(1+\left(\frac{\ol\Ph}{\ul\Ph}\right)^2\right) \right]\n{\Psi}{H^s(\Om)}.
\end{aligned}\end{equation}
The truncation operator $T:H^s(\Om)\ra H^s(\Om)$ satisfies

\begin{equation}\nn\begin{aligned}
\n{T}{Lip\left(H^s(\Om),H^s(\Om)\right)} = 1.
\end{aligned}\end{equation}
The proof is then complete. \cqfd


In the case of a contraction map, the iterative algorithm converges exponentially to the fixed point $a_\Psi$ and yields a map
\begin{equation}\nn\begin{aligned}
I:H^s(D) &\ra\A_s\\
\Psi &\mt a_\Psi.
\end{aligned}\end{equation}

From the Lipschitz continuity of $\ti G_\Psi$ with respect to $\Psi$, the following stability result holds. 
\begin{prop} For all $\Psi, \Psi^\prime \in H^s(D)$ such that $G_\Psi \circ F$ and $G_{\Psi^\prime} \circ F$ are contractions, we have 
$$
\| I[\Psi] - I[\Psi^\prime] \|_{H^s(D)} \leq C \| \Psi - \Psi^\prime\|_{H^s(D)}
$$ for some positive constant $C$. 
\end{prop}

\Proof Consider $\Psi, \Psi^\prime \in H^s(D)$ such that $G_\Psi \circ F$ and $G_{\Psi^\prime} \circ F$ are contractions and call $a_\Psi$ and $a_{\Psi^\prime}$ their fixed points. We have for any $\Ph\in\B_s$,

\begin{equation}\nn\begin{aligned}
\n{G_\Psi[\Ph]-G_{\Psi^\prime}[\Ph]}{H^s(\Om)}\leq \n{\ti G_\Psi[\Ph]-\ti G_{\Psi^\prime}[\Ph]}{H^s(\Om)} + \n{\frac{\Psi-\Psi^\prime}{\Ph^2}}{H^s(\Om)}.
\end{aligned}\end{equation}
Remarking that $u:=G_\Psi[\Ph]-G_{\Psi^\prime}[\Ph]$ satisfies

\begin{equation}\nn\left\{\begin{aligned}
\g\p(\Ph^2\g u) &=2\g\p\left[(\Psi-\Psi^\prime)\g\log\Ph\right]\quad \tin D,\\
u &=0\quad \ton\dr D,
\end{aligned}\right.\end{equation}
it follows that

\begin{equation}\nn\begin{aligned}
\n{\ti G_\Psi[\Ph]-\ti G_{\Psi^\prime}[\Ph]}{H^s(\Om)}\leq \frac {2R_{\B_s}} {\ul\Ph^3}\n{\Psi-\Psi^\prime}{L^2(\Om)}
\end{aligned}\end{equation}
and so

\begin{equation}\nn\begin{aligned}
\n{G_\Psi[\Ph]-G_{\Psi^\prime}[\Ph]}{H^s(\Om)}\leq \frac {4R_{\B_s}} {\ul\Ph^3}\n{\Psi-\Psi^\prime}{H^s(\Om)}.
\end{aligned}\end{equation}

We can now estimate

\begin{equation}\nn\begin{aligned}
\n{a_\Psi-a_{\Psi^\prime}}{H^s(\Om)} &=\n{G_\Psi \circ F[a_\Psi]-G_{\Psi^\prime} \circ F[a_{\Psi^\prime}]}{H^s(\Om)}\\
&\leq \n{G_\Psi \circ F[a_\Psi]-G_{\Psi^\prime} \circ F[a_{\Psi}]}{H^s(\Om)} + \n{G_{\Psi^\prime} \circ F[a_\Psi]-G_{\Psi^\prime} \circ F[a_{\Psi^\prime}]}{H^s(\Om)}\\
&\leq \n{G_\Psi-G_{\Psi^\prime}}{H^s(\Om)} + \n{G_{\Psi^\prime} \circ F[a_\Psi]-G_{\Psi^\prime} \circ F[a_{\Psi^\prime}]}{H^s(\Om)}\\
&\leq \frac {4R_{\B_s}} {\ul\Ph^3}\n{\Psi-\Psi^\prime}{H^s(\Om)}+\n{G_{\Psi^\prime} \circ F}{Lip\left(H^s(\Om)\right)}\n{a_\Psi-a_{\Psi^\prime}}{H^s(\Om)}.
\end{aligned}\end{equation}

Let $\kappa:=\n{G_{\Psi^\prime} \circ F}{Lip\left(H^s(\Om)\right)}<1$. It follows that

\begin{equation}\nn\begin{aligned}
\n{a_\Psi-a_{\Psi^\prime}}{H^s(\Om)} \leq \frac {4R_{\B_s}} {\ul\Ph^3(1-\kappa)}\n{\Psi-\Psi^\prime}{H^s(\Om)},
\end{aligned}\end{equation}
which completes the proof. \cqfd

\bigskip
\bigskip

\section{Numerical simulations}

In this section, we show how this new technique allows a very good reconstruction of highly discontinuous absorption map. We consider here a realistic absorption map taken from a blood vessels picture.

\subsection{Forward problem}
As we said in introduction, the main application of this acousto-optic method would be the imaging of red light absorption which has high contrast in tumors due to the high level of vascularization.

\begin{figure}[!h]
\begin{center}
\def\scale{0.78}
\def\size{5 cm}
\begin{tikzpicture}[scale= \scale]
\begin{axis}[width=1.6*\size, height=\size,axis on top, scale only axis, xmin=0, xmax=1.6, ymin=0, ymax=1]
\addplot graphics [xmin=0,xmax=1.6,ymin=0,ymax=1]  {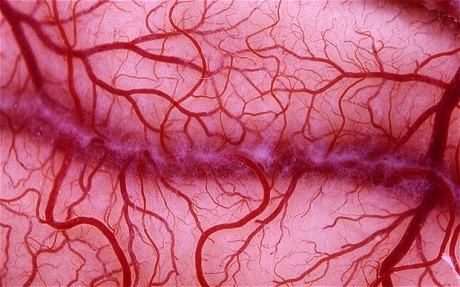};
\end{axis}
\end{tikzpicture}
\begin{tikzpicture}[scale= \scale]
\begin{axis}[width=1.6*\size, height=\size,axis on top, scale only axis, xmin=0, xmax=1.6, ymin=0, ymax=1, colormap/blackwhite, colorbar,point meta min=1,point meta max=1.98]
\addplot graphics [xmin=0,xmax=1.6,ymin=0,ymax=1]  {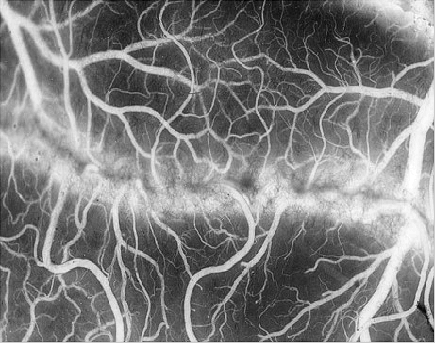};
\end{axis}
\end{tikzpicture}
\caption{\label{fig_aot41}Realistic biological light absorption map. (1) A real picture of living membrane by transparency. (2) The absorption map chosen for the numerical experiments. The resolution is about $132$k pixels.}
\end{center}
\end{figure}

\begin{figure}[!h]
\begin{center}
\def\scale{0.78}
\def\size{5 cm}
\begin{tikzpicture}[scale= \scale]
\begin{axis}[width=1.6*\size, height=\size,axis on top, scale only axis, xmin=0, xmax=1.6, ymin=0, ymax=1, colormap/blackwhite, colorbar,point meta min=1,point meta max=1.98]
\addplot graphics [xmin=0,xmax=1.6,ymin=0,ymax=1]  {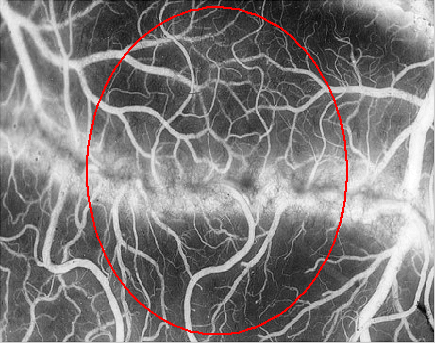};
\end{axis}
\end{tikzpicture}
\caption{\label{fig_aot42} Absorption map in $\Om$ and the domain of interest $D:=D\big((0.8,0.5),0.48\big)$ in red.}
\end{center}
\end{figure}

In the following, the domain is fixed to $\Om=]0,1.6[\times]0,1[$ and we consider the absorption map $a$ given by Figure \ref{fig_aot41} (2). We define our domain $D$ as a disk strictly included in $\Om$ represented by the red circle in Figure \ref{fig_aot42}.

\begin{figure}[!h]
\begin{center}
\def\scale{0.9}
\def\size{5 cm}
\begin{tikzpicture}[scale= \scale]
\begin{axis}[title = $M_\eta$, width=\size, height=\size,axis on top, scale only axis, xmin=0, xmax=128, ymin=0, ymax=0.8, colormap/jet, colorbar,point meta min=-0.01,point meta max=0.01,xlabel = $y$,ylabel = $r$]
\addplot graphics [xmin=0,xmax=128,ymin=0,ymax=0.8]  {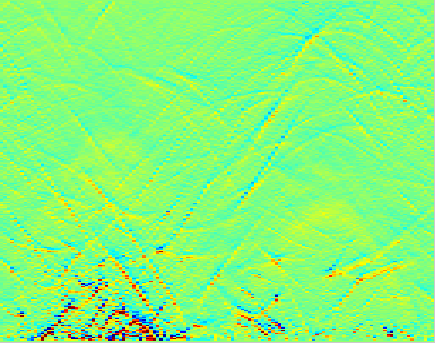};
\end{axis}
\end{tikzpicture}
\begin{tikzpicture}[scale= \scale]
\begin{axis}[title = $\Ra\Psi$, width=\size, height=\size, axis on top, scale only axis, xmin=0, xmax=128, ymin=0, ymax=0.8, colormap/jet, colorbar,point meta min=-0.00003,point meta max=0.00003,xlabel = $y$,ylabel = $r$]
\addplot graphics [xmin=0,xmax=128,ymin=0,ymax=0.8]  {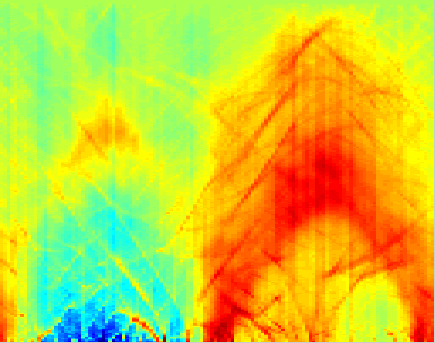};
\end{axis}
\end{tikzpicture}
\caption{\label{fig_aot43} Computed measurement $M_\eta(y,r)$ and the deduced approximation of $\Ra[\Psi]$. We used $128$ acoustic centers on $\dr D$.}
\end{center}
\end{figure}

Using the same method as in the numerical simulation in \cite{chap_2} we compute the forward problem in order to generate virtual measurements. For some centers $y$ taken on $Y:=\dr D$, $r>0$ and $\eta=10^{-4}$ fixed, we compute a discrete 
form of the map 

\begin{equation}\nn\begin{aligned}
v_{y,r,\eta}(x)=\frac \eta r w\left(\frac{|x-y|-r}\eta\right)\frac{x-y}{|x-y|},
\end{aligned}\end{equation}
where the wave shape $w$ is defined by

\begin{equation}\nn w(t)=\left\{\begin{aligned}
&\exp\left(\frac{1}{t^2-1}\right)\bb &t\in]-1,1[,\\
&0\bb &\text{otherwise}.
\end{aligned}\right.\end{equation}
From this map, we compute the displaced absorption as $a_v=a\circ(Id+v)^{-1}$ and the variation of the fluence $\Ph_v-\Ph$. Its cross correlation on the boundary leads to the measurement

\begin{equation}\nn\begin{aligned}
M_\eta(y,r)=\io(a_{v_{y,r,\eta}}-a)\Ph\Ph_{v_{y,r,\eta}},
\end{aligned}\end{equation}
represented in Figure \ref{fig_aot43} (1). From that, we apply Theorems \ref{theo_meas} and \ref{theo_step4} to get an approximation of $\Ra[\Psi]$ up to a function depending only on $y$. 

The non vertical visible lines on the illustration of $\Ra[\Psi]$ are due to the presences of blood vessels. The vertical lines are just numerical artifacts due to the integration. As we only need to know $\Ra[\Psi]$ up to a function depending only on $y$ to theoretically reconstruct the absorption, this last numerical issue is not important. Now, from numerical spherical means Radon transform inversion, we compute the internal data map $\Psi$ inside $D$.

\begin{figure}[!h]
\begin{center}
\def\scale{1}
\def\size{5 cm}
\begin{tikzpicture}[scale= \scale]
\begin{axis}[title = $\Psi$, width=\size, height=\size,axis on top, scale only axis, xmin=0.32, xmax=1.28, ymin=0.02, ymax=0.98, colormap/jet, colorbar,point meta min=-0.00019,point meta max=0.0001]
\addplot graphics [xmin=0.32, xmax=1.28, ymin=0.02, ymax=0.98]  {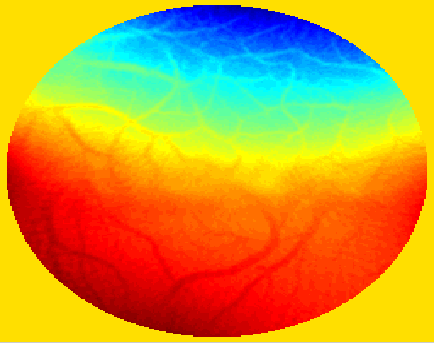};
\end{axis}
\end{tikzpicture}
\caption{\label{fig_aot44} Internal data map $\Psi$ computed inside the domain of interest $D$.}
\end{center}
\end{figure}

As we can observe the blood vessels in the representation of the map $\Psi$, we shall confirm that there is a good information about the absorption map. If we try the algorithm presented in \cite{chap_2}, we take the derivative of the data map $\Psi$ in order to compute the source term $\La\Psi$ which destroys the information due to the numerical noise. It is even worse with additional measurement noise. Here, we use the fixed point algorithm proposed in Theorem \ref{theo_contraction}. We compute the fixed point sequence $(a_n,\Ph_n)_{n\in\N}$ defined by

\begin{equation}\nn (a_0,\Ph_0):\left\{\begin{aligned}
a_0 &=1\bb\tin\Om,\\
\Ph_0 &:\left\{\begin{aligned}-\La\Ph_0+a_0\Ph_0 &=0\quad \tin\Om,\\
			      l\dr_\nu\Ph_0+\Ph_0 &=g\quad \ton\dr\Om,
	\end{aligned}\right.
\end{aligned}\right.\end{equation}
and
\begin{equation}\nn \forall \; n\in\N,\bb(a_{n+1},\Ph_{n+1}):\left\{\begin{aligned}
\ti a_{n+1} &:\left\{\begin{aligned}\g\p(\Ph_n^2\g\ti a_{n+1}) &=2\g\p(\Psi\g\log\Ph_{n})\quad \tin D,\\
			      \ti a_{n+1} &=-\frac{\Psi}{\Ph_{n}^2}\quad \ton\dr D,
	\end{aligned}\right.\\
a_{n+1} &:\left\{\begin{aligned} &1+\frac{\Psi}{\Ph_{n}^2}+\ti a_{n+1}\quad \tin D,\\
			      &1 \quad \tin\Om\bs\ol D,
	\end{aligned}\right.\\
\Ph_{n+1} &:\left\{\begin{aligned}-\La\Ph_{n+1}+a_{n}\Ph_{n+1} &=0\quad \tin\Om,\\
			      l\dr_\nu\Ph_{n+1}+\Ph_{n+1} &=g\quad \ton\dr\Om .
	\end{aligned}\right.
\end{aligned}\right.\end{equation}
After few iterations of this sequence, we get a good reconstruction of the absorption map $a$. To fix the ideas, let us say that the variational information about $a$ is in the map $\Psi/\Phi^2$. We correct it with a smooth function $\ti a$ in order to reach the map $a$. The difference of these two functions gives an approximation of $a-1$ at each iteration.

\begin{figure}[!h]
\begin{center}
\def\scale{1}
\def\size{5 cm}
\begin{tikzpicture}[scale= \scale]
\begin{axis}[title = $\Psi/\Phi^2$, width=\size, height=\size,axis on top, scale only axis, xmin=0.32, xmax=1.28, ymin=0.02, ymax=1.98, colormap/jet, colorbar,point meta min=-0.0025,point meta max= 0.0039]
\addplot graphics [xmin=0.32, xmax=1.28, ymin=0.02, ymax=1.98]  {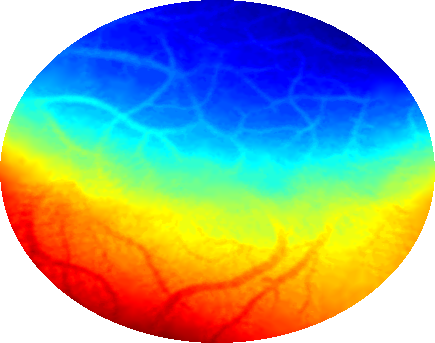};
\end{axis}
\end{tikzpicture}
\begin{tikzpicture}[scale= \scale]
\begin{axis}[title = $\ti a$, width=\size, height=\size,axis on top, scale only axis, xmin=0.32, xmax=1.28, ymin=0.02, ymax=0.98, colormap/jet, colorbar,point meta min= -0.0039,point meta max= 0.0025]
\addplot graphics [xmin=0.32, xmax=1.28, ymin=0.02, ymax=0.98]  {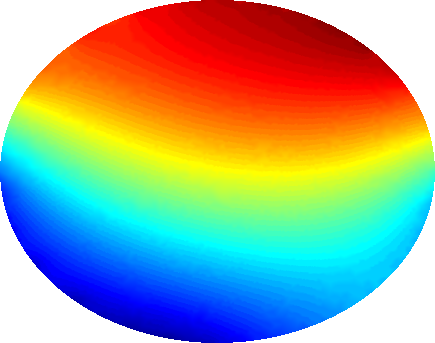};
\end{axis}
\end{tikzpicture}
\caption{\label{fig_aot44b} The map $\Psi/\Phi^2$ where the blood vessels are visible and the map $\ti a$ after $10$ iterations of the fixed point algorithm.}
\end{center}
\end{figure}

\begin{figure}[!h]
\begin{center}
\def\scale{1}
\def\size{5 cm}
\begin{tikzpicture}[scale= \scale]
\begin{axis}[title = (1), width=\size, height=\size,axis on top, scale only axis, xmin=0.32, xmax=1.28, ymin=0.02, ymax=1.98, colormap/jet, colorbar,point meta min=1,point meta max=1.78]
\addplot graphics [xmin=0.32, xmax=1.28, ymin=0.02, ymax=1.98]  {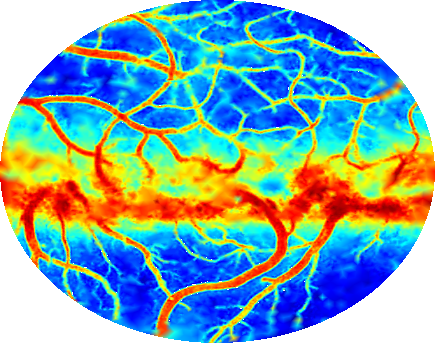};
\end{axis}
\end{tikzpicture}
\begin{tikzpicture}[scale= \scale]
\begin{axis}[title = (2), width=\size, height=\size,axis on top, scale only axis, xmin=0.32, xmax=1.28, ymin=0.02, ymax=0.98, colormap/jet, colorbar,point meta min=0.99,point meta max=1.58]
\addplot graphics [xmin=0.32, xmax=1.28, ymin=0.02, ymax=0.98]  {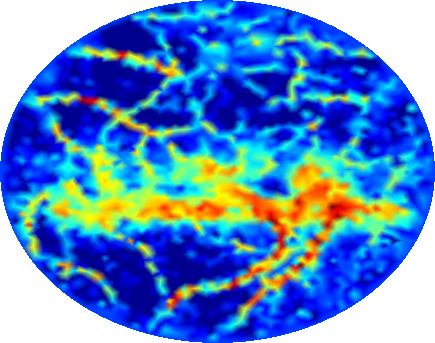};
\end{axis}
\end{tikzpicture}

\begin{tikzpicture}[scale= \scale]
\begin{axis}[title = (3), width=\size, height=\size,axis on top, scale only axis, xmin=0.32, xmax=1.28, ymin=0.02, ymax=0.98, colormap/jet, colorbar,point meta min=0.995,point meta max=1.67]
\addplot graphics [xmin=0.32, xmax=1.28, ymin=0.02, ymax=0.98]  {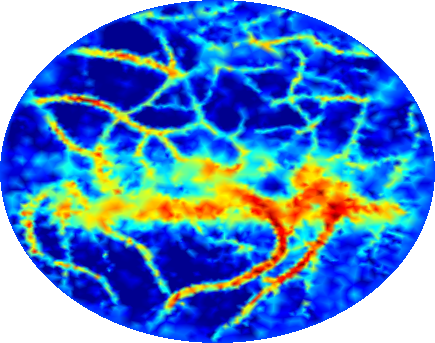};
\end{axis}
\end{tikzpicture}
\begin{tikzpicture}[scale= \scale]
\begin{axis}[title = (4),width=\size, height=\size,axis on top, scale only axis, xmin=0.32, xmax=1.28, ymin=0.02, ymax=0.98, colormap/jet, colorbar,point meta min=0.995,point meta max=1.72]
\addplot graphics [xmin=0.32, xmax=1.28, ymin=0.02, ymax=0.98]  {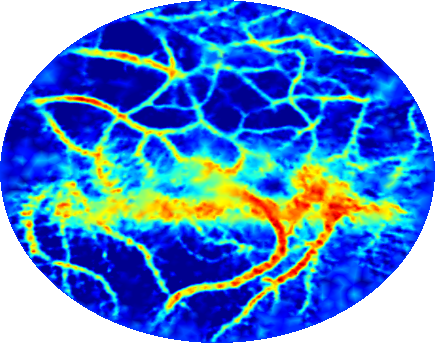};
\end{axis}
\end{tikzpicture}
\caption{\label{fig_aot45} Reconstruction of the absorption map after $10$ iterations of the fixed point sequence. (1) The true absorption. (2) Reconstruction using uniform mesh of $5$k triangles. (3) Reconstruction with non uniform mesh of $13$k triangles. (4) Reconstruction with non uniform mesh of $106$k triangles.}
\end{center}
\end{figure}

\begin{rem}
The power of this algorithm is that we avoid the derivation of the data map $\Psi$ and we only solve elliptic equation for smooth solutions $\Phi_n$ and $\ti a$. This provides a good reconstruction of the discontinuities of the absorption map $a$ and illustrates the fixed point Theorem \ref{theo_contraction} which works for functions in $H^s(\Om)$ with $s<1/2$.
\end{rem}

\begin{rem}
Our finest reconstruction is given in Figure \ref{fig_aot45} (4). Even if the vessels are easy to recognize, two problems occur. The first one is that the reconstructed solution is lightly attenuated. This is due to the approximation made using the asymptotic formula given in Theorem \ref{theo_step4}. A nice improvement would be to solve a deconvolution problem instead of the asymptotic formula. The second problem is the strong attenuation close to the boundary $\dr D$. This phenomenon is normal and is due to the fact that the measurements have no sense for small radius $r$. In the mathematical part, we have supposed that $a=a_0$ in a neighborhood of $\dr D$. In this numerical example, this hypothesis is not respected and the consequence is that the reconstruction is not valid close to $\dr D$. Nevertheless, the inside part of the reconstruction is quite satisfying. 
\end{rem}

\section{Concluding remarks}
In this paper we have introduced for the first time a mathematical and numerical framework for reconstructing highly discontinuous contrast distributions from internal measurements. 
The framework yields stable and accurate reconstructions. We have illustrated our approach on a 
highly discontinuous absorption map, chosen from a real biological tissue data. Many challenging problems are still open. It would be very interesting to develop an optimal control scheme for reconstructing highly discontinuous contrast distributions and prove its convergence, starting from a good initial guess. Another challenging problem is to estimate the resolution of the developed approach in terms of the signal-to-noise ratio in the data.

\appendix

\section{Spherical density of $Da$}\label{app_density} 

\begin{lem} Consider $a\in SBV^\ii(\Om)$ constant out of the convex $<D\subset \Om$ and the mollifier sequence $w_\eta(r)=\frac 1 \eta w\left(\frac 1\eta\right)$. Suppose that $Y$ satisfies the wrap condition around $D$, then the sequence of functions defined on $\Sigma$

\begin{equation}\nn\begin{aligned}
\ph_\eta(y,r)=\int_\Om w_\eta(|x-y|-r)|D a|(dx)
\end{aligned}\end{equation}
satisfies

\begin{equation}\nn\begin{aligned}
\n{\ph_\eta}{L^2(\Sigma)}\leq C\eta^{-\frac 1{2d}},
\end{aligned}\end{equation}
where $C$ depends on $|D a|(\Om)$, $|Y|$, and the wrap constant. 
\end{lem}

\Proof We develop $\n{\ph_\eta}{L^2(\Sigma)}^2$ norm as
\begin{equation}\nn\begin{aligned}
\n {\ph_\eta}{L^2(\Si)}^2 &=\int_\Si\int_\Om\int_\Om w_\eta(|x-y|-r)w_\eta(|x'-y|-r)| Da|(dx)| Da|(dx')dydr\\
&=\int_Y\int_\Om\int_\Om \int_0^R w_\eta(|x-y|-r)w_\eta(|x'-y|-r)dr| Da|(dx)| Da|(dx')dy\\
&=\int_Y\int_\Om\int_\Om\ol w_\eta(|x-y|-|x'-y|)| Da|(dx)| Da|(dx')dy, \\
\end{aligned}\end{equation}
where

\begin{equation}\nn\begin{aligned}
\ol w_\eta(r)=\int_\R w_\eta(r-\rho)w_\eta(-\rho)d\rho
\end{aligned}\end{equation}
satisfies $\supp(\ol w_\eta)\su[-2\eta,2\eta]$, $\n{\ol w_\eta}{L^1(\R)}\leq 1$ and $\ol w_\eta\leq\frac 1 \eta$.
Let us fix $\e>0$ and define $Z_\e=\{(x,x')\in\Om^2,\ |x-x'|\leq\e\}$. First, we have

\begin{equation}\nn\begin{aligned}
\int_Y\int_{Z_\e}\ol w_\eta(|x-y|-|x'-y|)| Da|(dx)| Da|(dx')dy &\leq \frac 1 \eta|Y|\int_\Om\int_{B(x,\e)}| Da|(dx')| Da|(dx)\\
&\leq \frac 1 \eta|Y|\int_\Om| Da|(B(x,\e))| Da|(dx).\\
\end{aligned}\end{equation}
Using the fact that $a\in SBV^\ii(\Om)$, the Radon measure $| Da|$ can be decomposed as
\begin{equation}\nn\begin{aligned}
| Da|=|\g_la|{\cal L}^d+|[a]_S|\cH_S^{d-1},
\end{aligned}\end{equation}
where $|\g_la|\in L^\ii(\Om)$ and $|[a]_S|\in L^\ii(S)$. Thus, we can control the upper (d-1)-densities of $| Da|$ using that for any $x\in\Om'$

\begin{equation}\nn\begin{aligned}
\frac 1 {\e^{d-1}}| Da|(B(x,\e))\leq \n {\g_la}{L^\ii(\Om)}\om^d\e+\n {[A]_S}{L^\ii(S)}\frac 1 {\e^{d-1}}\cH^{d-1}(S\cap B(x,\e)).
\end{aligned}\end{equation}
In fact, \cite[Theorem 6.2]{mattila} says that for any $x\in S$,

\begin{equation}\nn\begin{aligned}
\limsup_{\e\to 0} \frac 1 {\e^{d-1}}\cH^{d-1}(S\cap B(x,\e)) &\leq 2^{d-1}\ \ &a.e. \ton S,\\
\limsup_{\e\to 0} \frac 1 {\e^{d-1}}\cH^{d-1}(S\cap B(x,\e)) &= 0\ \ &a.e. \ton \Om'\bs S,
\end{aligned}\end{equation}
which implies for $| Da|$ that

\begin{equation}\nn\begin{aligned}
\limsup_{\e\to 0} \frac 1 {\e^{d-1}}| Da|(B(x,\e)) &\leq \n {[A]_S}{L^\ii(S)}2^{d-1}\ \ &a.e. \ton S,\\
\limsup_{\e\to 0} \frac 1 {\e^{d-1}}| Da|(B(x,\e)) &= 0\ \ &a.e. \ton \Om'\bs S.
\end{aligned}\end{equation}
Using Fatou lemma, it follows that
\begin{equation}\nn\begin{aligned}
\limsup_{\e\to 0}\int_{\Om'}\frac 1 {\e^{d-1}}| Da|(B(x,\e))| Da|(dx) &\leq \int_{\Om'}\limsup_{\e\to 0}\frac 1 {\e^{d-1}}| Da|(B(x,\e))| Da|(dx)\\
&\leq \n {[A]_S}{L^\ii(S)}2^{d-1}\cH^{d-1}(S).
\end{aligned}\end{equation}
That simply shows that the left-hand integral is bounded when $\e$ goes to zero. We finally arrive at

\begin{equation}\label{eq_proof1}\begin{aligned} 
\int_Y\int_{Z_\e}\ol w_\eta(|x-y|-|x'-y|)| Da|(dx)| Da|(dx')dy &\leq C_1\frac {\e^{d-1}} \eta, 
\end{aligned}\end{equation}
where the constant $C_1$ depends on $| Da|(\Om)$ and $|Y|$. The second integral that we have to control is
\begin{equation}\nn\begin{aligned}
\int_Y\int_{\Om'^2\bs Z_\e}\ol w_\eta(|x-y|-|x'-y|)| Da|(dx)| Da|(dx')dy.
\end{aligned}\end{equation}
For that, we define for any $(x,x')\in\Om^2$ the set
\begin{equation}\nn\begin{aligned}
Y_\eta(x,x')=\left\{y\in Y,\ \left||x-y|-|x'-y|\right|\leq 2\eta\right\}.
\end{aligned}\end{equation}
As $Y$ satisfies the wrap condition around $\Om'$, a computation leads to 

\begin{equation}\nn\begin{aligned}
\cH^{d-1}\left(Y_\eta(x,x')\right)\leq 2C_2\frac \eta {\e}\ \ \ \ \forall \; (x,x')\in\Om'^2\bs Z_\e,
\end{aligned}\end{equation}
where $C_2$ is the wrap constant relative to $Y$ and $\Om'$. We can now control the second term, 
\begin{equation}\label{eq_proof2}\begin{aligned}
\int_Y\int_{\Om'^2\bs Z_\e}\ol w_\eta(|x-y|-|x'-y|)| Da|(dx)| Da|(dx')dy &\leq \frac 1\eta \int_{\Om'^2}\cH^{d-1}\left(Y_\eta(x,x')\right)| Da|(dx)| Da|(dx')\\
&\leq \frac{C_2} \e| Da|(\Om)^2.
\end{aligned}\end{equation}
Finally, putting together (\ref{eq_proof1}) and (\ref{eq_proof2}), we obtain
\begin{equation}\nn\begin{aligned}
\n {\ph_\eta}{L^2(\Si)}^2 &\leq C_1\frac {\e^{d-1}} \eta+2\frac{C_2} \e| Da|(\Om)^2,
\end{aligned}\end{equation}
which is true for any choice of $\e>0$. So, we fix it at the best choice $\e=\eta^{1/d}$ to obtain
\begin{equation}\nn\begin{aligned}
\n {\ph_\eta}{L^2(\Si)}^2 &\leq (C_1+2C_2| Da|(\Om)^2)\eta^{-\frac 1 d},
\end{aligned}\end{equation}
which concludes the proof.\cqfd
\bigskip

\section{Sobolev spaces with fractional order and Helmholtz decomposition} \label{appendsobolev}

On the smooth open domain $D$ of $\R^d$, for any $\al\geq 0$ the Sobolev space $H^\al(D)$ is defined as usual. We shall also consider the space of functions of $H^\al(D)$ supported in a compact $K$ denoted $H^\al_K(D)$. As the functions of $H^\al_K(D)$ can be extended by zero outside of $D$, we can define their Fourier transform and use the following characterization,

\begin{de}For any $\al\geq 0$, $K\su D$ compact we define
\begin{equation}\nn\begin{aligned}
H^\al_K(D)=\left\{f\in L^2(D),\ \supp(u)\su K,\ \int_{\R^d}|\wh f|^2(\xi)(1+|\xi|^2)^\al d\xi<+\ii\right\}
\end{aligned}\end{equation}
and for any $f\in H^\al_K(D)$ we will denote $$\n{f}{H^\al(D)}=\left(\frac 1 {(2\pi)^d}\int_{\R^d}|\wh f|^2(\xi)(1+|\xi|^2)^\al d\xi\right)^{\frac 1 2}.$$
\end{de}
We define now $H_K^{-\al}(D)$ by duality.

\begin{de}For any $\al> 0$, $K\su D$ compact we define
\begin{equation}\nn\begin{aligned}
H^{-\al}_K(D)=\left\{u\in H^\al(D)',\ \supp(u)\su K\right\}
\end{aligned}\end{equation}
endowed with the continuity norm.
\end{de}
Fortunately, these spaces have also a Fourier characterization. For any $u\in H^{-\al}_K(D)$, $u$ is a compact supported distribution, i.e., an element of $\E'(D)$, which naturally embeds in $\Sa'(\R^d)$. So, the Fourier transform $\wh u$ is defined in $\Sa'(\R^d)$.

\begin{prop}For any $\al> 0$, $K\su D$ compact,
\begin{equation}\nn\begin{aligned}
H^{-\al}_K(D)=\left\{u\in \E'(D),\ \supp(u)\su K,\ \wh u\in L^1_{\mathrm{loc}}\big(\R^d\big),\ \int_{\R^d}|\wh u|^2(\xi)(1+|\xi|^2)^{-\al} d\xi<+\ii\right\}.
\end{aligned}\end{equation}
\end{prop}
\Proof Let us take $u\in H^{-\al}_K(D)$. As $u\in\Sa'(R^d)$, we take $\wh u\in\Sa'(\R^d)$, $\ph\in\Sa'(\R^d)$ and we compute
\begin{equation}\nn\begin{aligned}
\left|\left<(1+|\xi|^2)^{-\al/2}\wh u,\ph\right>_{\Sa'(\R^d),\Sa(\R^d)} \right|&=\left|\left<u,
\widehat{\left[{(1+|x|^2)^{-\al/2}\ph}\right]}\right>_{\Sa'(\R^d),\Sa(\R^d)}\right|\\
&\leq \n u{H^\al(D)'}\n{\widehat{\left[{(1+|x|^2)^{-\al/2}\ph}\right]}}{H^\al(D)}\\
&\leq (2\pi)^{d/2}\n u{H^\al(D)'}\n{\ph}{L^2(D)},\\
\end{aligned}\end{equation}
which proves that $(1+|\xi|^2)^{-\al/2}\wh u\in L^2(\R^d)$ and $$\left(\frac 1{(2\pi)^d}\int_{\R^d}|\wh u|^2(\xi)(1+|\xi|^2)^{-\al} d\xi\right)^{1/2}\leq \n u{H^\al(D)'}.$$ Conversely,  if $u$ satisfies these conditions, we show that it is in $H^\al(D)'$ and that $$\n u{H^\al(D)'}\leq\left(\frac 1{(2\pi)^d}\int_{\R^d}|\wh u|^2(\xi)(1+|\xi|^2)^{-\al} d\xi\right)^{1/2}.$$
Then the proof is complete. \cqfd

We can now define the Helmholtz decomposition of a distribution vectorial field in the Sobolev sense for fractional order greater than $-1$. This allows us to precise the regularity of $\Psi$ depending on the regularity of $a$.
\bigskip

\section{Kernel operators in partial Sobolev spaces}

In this appendix, we give two useful results about some kernel operators acting on one variable of a function. These results are given for functions defined in $\R^d$ in order to use the Fourier transform. They stay valid for functions defined on any manifold isomorphic to an open domain of $\R^d$ up to a multiplicative constant depending on the isomorphism.

\begin{lem}\label{lem_ker1} Consider a kernel $\theta\in L^2(\R^2)$ and the operator $T:L^2(\R^d)\ra L^2(\R^d)$ defined by
\begin{equation}\nn\begin{aligned}
T[f](x)=\int_\R f(t,\ti x)\theta(t,x_1)dt
\end{aligned}\end{equation}
for a.e. $x\in\R^d$ with $\ti x=(x_2,\dots,x_d)$. If, for $s>0$, $f\in H^s(\R^d)$ and $\theta\in H^s(\R^2)$, then $T[f]\in H^s(\R^d)$ and we have
\begin{equation}\nn\begin{aligned}
\n{T[f]}{H^s(\R^d)}\leq\n{\theta}{H^s(\R^2)}\n{f}{H^s(\R^d)}.
\end{aligned}\end{equation}
\end{lem}

\Proof Let us compute the Fourier transform of $T[f]$,
\begin{equation}\nn\begin{aligned}
\wh{T[f]}(\xi) &= \int_\R\int_\R\int_{\R^{d-1}}f(t,\ti x)\theta(t,x_1)e^{-i x_1\xi_1}e^{-i \ti x\cdot\ti \xi}d\ti x d x_1dt\\
&= \int_\R\overset{\ \ti x}{\wh{f}}(t,\ti\xi)\overset{\ x_1}{\wh\theta}(t,\xi_1)dt,
\end{aligned}\end{equation}
so
\begin{equation}\nn\begin{aligned}
|\wh{T[f]}|^2(\xi) &\leq\int_\R |\overset{\ \ti x}{\wh{f}}(t,\ti\xi)|^2dt\int_\R|\overset{\ x_1}{\wh\theta}(t,\xi_1)|^2dt.
\end{aligned}\end{equation}
Then, using Plancherel theorem,
\begin{equation}\nn\begin{aligned}
\int_\R |\overset{\ \ti x}{\wh{f}}(t,\ti\xi)|^2dt=\frac 1{2\pi}\int_{\R}|\wh f(\xi)|^2d\xi_1
\end{aligned}\end{equation}
and
\begin{equation}\nn\begin{aligned}
\int_\R|\overset{\ x_1}{\wh\theta}(t,\xi_1)|^2dt= \frac 1{2\pi}\int_{\R}|\wh \theta(\tau,\xi_1)|^2d\tau.
\end{aligned}\end{equation}
Hence,
\begin{equation}\nn\begin{aligned}
|\wh{T[f]}|^2(\xi) &\leq \frac 1{(2\pi)^2}\int_{\R}|\wh f(\xi)|^2d\xi_1\int_{\R}|\wh \theta(\tau,\xi_1)|^2d\tau\\
|\wh{T[f]}|^2(\xi)\left(1+|\xi|^2\right)^s &\leq \frac 1{(2\pi)^2}\int_{\R}|\wh f(\xi)|^2\left(1+|\ti\xi|^2\right)^sd\xi_1\int_{\R}|\wh \theta(\tau,\xi_1)|^2\left(1+\xi_1^2\right)^sd\tau\\
\frac1{(2\pi)^d}\int_{\R^d}|\wh{T[f]}|^2(\xi)\left(1+|\xi|^2\right)^sd\xi &\leq\\ \frac 1{(2\pi)^d} \int_{\R^d}|\wh f(\xi)|^2\left(1+|\xi|^2\right)^s &d\xi\frac1{(2\pi)^2} \int_{\R^2}|\wh \theta(\tau,\xi_1)|^2\left(1+\xi_1^2+\tau^2\right)^sd\tau d\xi_1,
\end{aligned}\end{equation}
which completes the proof.  \cqfd

In the case where the kernel is approaching a delta function, it is useful to understand how the operator is approaching the identity.

\begin{lem}\label{lem_ker2} Consider $w\in\C^\ii_c(\R)$ supported in $[-1,1]$, non negative and satisfying $\n w{L^1(\R)}=1$. For any $\eta>0,\ t\in\R$ we denote $w_\eta(t)=\frac 1\eta w\left(\frac t \eta\right)$. Let us consider the sequence of operator $T_\eta:L^2(\R^d)\ra L^2(\R^d)$ defined by
\begin{equation}\nn\begin{aligned}
T_\eta[f](x)=\int_\R f(t,\ti x)w_\eta(x_1-t)dt.
\end{aligned}\end{equation}

For all $\al\geq 0$ and $\eta>0$, $T_\eta$ is continuous operator $T_\eta:H^\al(\R^d)\ra H^\al(\R^d)$ and for all $\beta>0$, $f\in H^{\al+\beta}(\R^d)$, $T_\eta[f]$ converges to $f$ in $H^\al(\R^d)$. More precisely,
\begin{equation}\nn\begin{aligned}
\n{T_\eta[f]-f}{(H^{\al}(\R^d)}\leq 2\eta^\frac{\be}{\be+1}\n{f}{H^{\al+\be}(\R^d)}.
\end{aligned}\end{equation}
\end{lem}

\Proof Let us compute the Fourier transform of $T[f]$,
\begin{equation}\nn\begin{aligned}
\wh{T[f]}(\xi) &= \int_\R\int_\R\int_{\R^{d-1}}f(t,\ti x)w_\eta(x_1-t)e^{-i x_1\xi_1}e^{-i \ti x\cdot\ti \xi}d\ti x d x_1dt\\
&= \int_\R\int_\R\int_{\R^{d-1}}f(t,\ti x)w_\eta(u)e^{-i u\xi_1}e^{-i t\xi_1}e^{-i \ti x\cdot\ti \xi}d\ti x d udt\\
&=\wh f(\xi)\wh{w_\eta}(\xi_1),
\end{aligned}\end{equation}
where $\wh{w_\eta}\leq 1$. This proves that $\n{T_\eta[f]}{H^\al(\R^d)}\leq\n{f}{H^\al(\R^d)}$. Now consider $\be>0$ and $f\in H^{\al+\be}(\R^d)$, we have
\begin{equation}\nn\begin{aligned}
\left(\wh{T[f]}-\wh f\right)(\xi) &=f(\xi)\int_\R w_\eta(t)(e^{-it\xi1}-1)dt,\\
\left|\widehat{T_\eta[f]}-\widehat{f}\right|^2(\xi) &\leq|\widehat{f}|^2(\xi)\int_\R w_\eta(t)|e^{-it\xi_1}-1|^2dt
\end{aligned}\end{equation}
by convexity, and we write,
\begin{equation}\nn\begin{aligned}
\left|\widehat{T_\eta[f]}-\widehat{f}\right|^2(\xi) &=|\widehat{f}|^2(\xi)\sup_{|t|\leq\eta}|e^{-it\xi_1}-1|^2.\\
\end{aligned}\end{equation}
A study of the function $\xi_1\mapsto\sup_{|z|\leq\eta}|e^{-iz\xi_1}-1|^2$ gives us that
\begin{equation}\nn\begin{aligned}
\sup_{|t|\leq\eta}|e^{-it\xi_1}-1|^2\leq 4\eta^{\frac{2\be}{\be+1}}\left(1+|\xi|^2\right)^\be, 
\end{aligned}\end{equation}
and we finally get

\begin{equation}\nn\begin{aligned}
\int_{\R^d}\left|\widehat{T_\eta[f]}-\widehat{f}\right|^2(\xi)\left(1+|\xi|^2\right)^\al d\xi\leq 4\eta^{\frac{2\be}{\be+1}}\int_{\R^d}|\widehat{f}|^2(\xi)\left(1+|\xi|^2\right)^{\al+\be} d\xi,
\end{aligned}\end{equation}
which is equivalent to
\begin{equation}\nn\begin{aligned}
\n{T_\eta[f]-f}{H^\al(\R^d)}\leq 2\eta^{\frac{\be}{\be+1}}\n{f}{H^{\al+\be}(\R^d)}.
\end{aligned}\end{equation}
Hence, the proof is complete. \cqfd


\begin{thebibliography}{99}

\bibitem{alberti_amb} G. Alberti and C. Mantegazza. A note on the theory of SBV functions. {\em Boll. Un. Mat. Ital.}, B (7)  11  (1997),  no. 2, 375--382.


\bibitem{AMMARI-08}
{ H.~Ammari}. {\em An Introduction to Mathematics of Emerging
Biomedical Imaging}. Vol. {62}, { Mathematics and Applications},
Springer-Verlag, Berlin, 2008.



\bibitem{chap_2}
{ H.~Ammari, E. Bossy, J. Garnier, L. H. Nguyen and L.
Seppecher}. {A
  reconstruction algorithm for ultrasound-modulated diffuse optical
  tomography}. {\em Proc. Amer. Math. Soc.}, to appear.


\bibitem{chap_1}
{ H.~Ammari, E. Bossy, J. Garnier, and L. Seppecher}. {
  Acousto-electromagnetic tomography}. {\em SIAM J. Appl. Math.}, 72
  (2012), 1592--1617.
  

\bibitem{chap_3} H. Ammari, J. Garnier, L.H. Nguyen, and L. Seppecher. Reconstruction of a piecewise smooth
absorption coefficient by an acousto-optic process. {\em Comm. Part. Differ. Equat.}, 
  38  (2013),  no. 10, 1737--1762. 


\bibitem{simon} { S.R. Arridge}. {Optical tomography in medical
imaging}. {\em Inverse Problems}, 15 (1999),  R41--R93.


\bibitem{born} {M. Born and E. Wolf}. {\em Principles of Optics}.  Cambridge University
Press, Cambridge, 1999.




%

\bibitem{fink} M. Fink and M. Tanter. Multiwave imaging and super resolution. {\em Phys. Today}, 63 (2010), 28--33. 



\bibitem{GilbargTrudinger:1977}
{ D.~Gilbarg and N.~S. Trudinger}. {\em Elliptic partial differential
  equations of second order}. Springer-Verlag, Berlin,
  1977.






%



\bibitem{mattila} P. Mattila. {\em Geometry of Sets and Measures in Euclidian Spaces. Fractals and rectifiability}. Cambridge Press, 1995. 


\bibitem{otmar} W. Naetar and O. Scherzer. Quantitative photoacoustic tomography with piecewise constant material parameters. Arxiv: 1403.2620. 

 

\bibitem{victoire} V. Palamodov. Remarks on the general Funk transform and thermoacoustic tomography. {\em Inverse Probl. Imaging}, 4  (2010),  no. 4, 693--702.

%




\bibitem{quinto} E.T. Quinto. Support theorems for the spherical
Radon transform on manifolds. {\em Int. Math. Res. Lett.}, 2006,
1--17 (Article ID 67205).

\bibitem{john} {J. C. Schotland}. {Direct reconstruction
methods in optical tomography}. {\em Lecture Notes in Math.}, Vol. 2035,
1--29, Springer-Verlag, Berlin, 2011.

\bibitem{seobook} J.K. Seo and E.J. Woo. \textsl{Nonlinear Inverse Problems in Imaging}. Wiley, 2013.

\bibitem{helmdec} W. Spr\"{o}ssig. On Helmholtz decompositions and their generalizations—an overview, Math. Methods Appl. Sci.,  33  (2010), 374--383.


\end{thebibliography}
\end{document}